\documentclass[11pt]{article}

\usepackage{amssymb}
\usepackage{amsmath}
\usepackage{theorem}
\usepackage{epsfig}
\usepackage{verbatim}
\usepackage{graphicx}

\usepackage{color}

\newtheorem{thm}{Theorem}[section]
\newtheorem{cor}[thm]{Corollary}

\newtheorem{rem}[thm]{Remark}

\newtheorem{lemma}[thm]{Lemma}

\newcommand{\R}{\mathbb{R}}
\newcommand{\C}{\mathbb{C}}
\newcommand{\Z}{\mathbb{Z}}

\newcommand{\T}{\mathbb{T}}
\newcommand{\D}{\displaystyle}

\newcommand{\grad}{\nabla}

\newcommand{\gradp}{\grad^{\bot}}

\newcommand{\da}{\partial_\alpha}

\newcommand{\la}{\Lambda}

\newcommand{\al}{\alpha}
\newcommand{\ep}{\varepsilon}

\newcommand{\ff}{\mathbb}
\newcommand{\pa}{\partial}

\newcommand{\be}{\beta}
\newcommand{\g}{\gamma}
\newcommand{\ze}{\zeta}

\begin{document}

\title{Rayleigh-Taylor breakdown for the Muskat problem\\
 with applications to water waves}
\date{\today}
\author{\'Angel Castro, Diego C\'ordoba, Charles Fefferman,\\
 Francisco Gancedo and Mar\'ia L\'opez-Fern\'andez.}
\maketitle

\begin{abstract}
The Muskat problem models the evolution of the interface between
two different fluids in porous media. The Rayleigh-Taylor
condition is natural to reach  linear stability of the Muskat
problem. We show that the Rayleigh-Taylor condition may hold
initially but break down in finite time. As a consequence of the
method used, we prove the existence of water waves turning.
\end{abstract}

\maketitle


\section{Introduction}

The Muskat problem \cite{Muskat} models the evolution of an
interface between two fluids of different characteristics in porous
media by means of Darcy's law:
\begin{equation}\label{darcy}
\frac{\mu}{\kappa}u=-\nabla p-(0,\mathrm{g}\rho),
\end{equation}
where $(x,t)\in\mathbb{R}^2\times\R^+$, $u = (u_1(x,t), u_2(x,t))$
is the incompressible velocity (i.e. $\nabla\cdot u =0$), $p=p(x,t)$
is the pressure, $\mu(x,t)$ is the dynamic viscosity, $\kappa$ is
the permeability of the isotropic medium, $\rho=\rho(x,t)$ is the
liquid density, and $\mathrm{g}$ is the acceleration due to gravity.
More precisely, the interface separates the domains $\Omega^1$ and
$\Omega^2$ defined by
\begin{equation*}
(\mu,\rho)(x_1,x_2,t)=\left\{\begin{array}{cl}
                    (\mu^1, \rho^1),& x\in\Omega^1(t)\\
                    (\mu^2, \rho^2),& x\in\Omega^2(t)=\mathbb{R}^2 - \Omega^1(t),
                 \end{array}\right.
\end{equation*} and $\mu^1, \mu^2, \rho^1, \rho^2$ are constants.
This physical situation is also related to the evolution of two fluids of
different characteristics in a Hele-Shaw cell \cite{H-S}, due to the
fact that the laws which model both phenomena are mathematically
analogous \cite{S-T}.

This paper is concerned with the case $\mu^1=\mu^2$ which provides
weak solutions of the following transport equation
\begin{align}
\begin{split}\label{trasporte}
\rho_t+u\cdot\grad\rho=0, \\
\rho_0 = \rho(x,0),\quad\quad x\in \mathbb{R}^2,
\end{split}
\end{align}
where initially the scalar $\rho_0$ is given by
\begin{equation}\label{frho}
 \rho_0= \rho(x_1,x_2,0)=\left\{\begin{array}{cl}
                    \rho^1\quad\mbox{in}&\Omega^1(0)=\{x_2>f_0(x_1)\}\\
                    \rho^2\quad\mbox{in}&\Omega^2(0)=\{x_2<f_0(x_1)\}.
                 \end{array}\right.
\end{equation}

Let the free boundary be parametrized by
$$\partial\Omega^j (t) =
\{z(\al, t) = (z_1(\al,t), z_2(\al,t)) : \al \in \R\}
$$
where
\begin{equation*}\label{curve}
z(\al,t)-(\al,0)
\end{equation*}
is $2\pi$-periodic in the space parameter $\alpha$ or, an open contour
vanishing at infinity
$$\D\lim_{\al\rightarrow\pm\infty}(z(\al,t)-(\al,0))=0$$
with initial data $z(\al, 0) = z_0(\al)=(\al,f_0(\al))$.
 From Darcy's law, we find that the vorticity is concentrated on the free boundary
$z(\al,t)$, and is given by a Dirac distribution as follows:
$$\gradp\cdot u(x,t)=\omega(\al,t)\delta(x-z(\al,t)),$$ with $\omega(\al,t)$ representing the vorticity strength i.e.
$\gradp\cdot u$ is a measure defined by
\begin{equation*}
<\gradp\cdot u,\eta>=\int \omega(\al,t)\eta(z(\al,t)) d\al,
\end{equation*}
with $\eta(x)$ a test function.

Then $z(\al, t)$ evolves with an incompressible velocity field
coming from the Biot-Savart law:
\begin{equation*}\label{BS}
u(x,t)=\gradp\Delta^{-1}\gradp\cdot u(x,t).
\end{equation*}

As $(x,t)$ approaches a  point  $z(\al, t)$ on the contour the velocity $u$ agrees, modulo tangential terms, with the  Birkhoff-Rott integral:
\begin{equation*}\label{BR}
BR(z,\omega)(\al,t)=\frac{1}{2\pi}PV\int
\frac{(z(\al,t)-z(\beta,t))^\bot}{|z(\al,t)-z(\beta,t)|^2}\omega(\beta,t)
d\beta.
\end{equation*}

This yields an appropriate contour dynamics system:
  \begin{align}
\begin{split}\label{fullpm}
z_t(\al,t)&=BR(z,\omega)(\al,t)+c(\al,t)\da z(\al,t),
\end{split}
\end{align}
where the term c represents the change of parametrization and does
not modify the geometric evolution of the curve \cite{HLS}.

 The well-posedness
is not guaranteed in general, in fact such a result turns out to be
false for some initial data. Rayleigh \cite{Ray} and Saffman-Taylor
\cite{S-T} gave a condition that must be satisfied for the
linearized model in order to have a solution locally in time, namely that the
normal component of the pressure gradient jump at the interface has
to have a distinguished sign. This is known as the Rayleigh-Taylor
condition:
$$\sigma(\al,t) = -(\nabla p^2(z(\al,t),t) - \nabla
p^1(z(\al,t),t))\cdot\partial^{\perp}_{\al}z(\al,t)>0,$$
where $\nabla p^j(z(\al,t),t)$ denotes the limit gradient of the pressure obtained approaching the boundary in the normal direction inside $\Omega^j(t)$.
We call $\sigma(\al,t)$ the Rayleigh-Taylor of the solution $z(\alpha,t)$.

Understanding the problem as weak solutions of
(\ref{darcy}-\ref{trasporte}) plus the incompressibility of the
velocity, we find that the continuity of the pressure
($p^2(z(\al,t),t)=p^1(z(\al,t),t)$) follows as a mathematical
consequence, making unnecessary to impose it as a physical
assumption (for more details see \cite{DY} and \cite{ADY}). For
the surface tension case, there is a jump discontinuity  of the
pressure across the interface which is modeled to be equal to the
local curvature times the surface tension coefficient:
$$
p^2(z(\al,t),t)-p^1(z(\al,t),t))=\tau \kappa(\al,t).
$$
This is known as the Laplace-Young condition, which makes the
initial value problem more regular. Then there are no instabilities
\cite{Esch1} but fingering phenomena arise \cite{Otto,Esch2}.

By means of Darcy's law, we can find the following formula for the
difference of the gradients of the pressure in the normal direction
and the strength of the vorticity:
$$\sigma(\al,t)=(\rho^2-\rho^1)\da z_1(\al,t)$$
 \begin{equation}\label{vae}
\omega(\al,t)=- (\rho^2-\rho^1)\da z_2(\al,t).
\end{equation}
Above $\mathrm{g}$ is taken equal to 1 for the sake of simplicity.

 Then, if we  choose an appropriate term c in equation (\ref{fullpm}) (see section
\ref{csle} below), the dynamics of the interface satisfies
\begin{equation}\label{ec1d}
z_t(\al,t) =\frac{\rho^2-\rho^1}{2\pi}PV\int
\frac{(z_1(\al,t)-z_1(\beta,t))}{|z(\al,t)-z(\beta,t)|^2}(\da
z(\al,t) - \da z(\beta,t)) d\beta.
\end{equation}

A wise choice of parametrization of the curve is to have $\da
z_1(\al,t)=1$ (for more details see \cite{DY}). This yields the
denser fluid below the less dense fluid if $\rho^2>\rho^1$ and
therefore the Rayleigh-Taylor condition holds as long as the
interface is a graph. This fact has been used in \cite{DY} to show
local existence in the stable case ($\rho^2 > \rho^1$), together
with ill-posedness in the unstable situation ($\rho^2 < \rho^1$).
Local existence for the general case ($\mu^1\neq\mu^2$) is shown
in \cite{ADY}, which was also treated in \cite{Yi,Am}.

>From \eqref{ec1d} it is easy to find the evolution equation for the
graph:
\begin{align}
\begin{split}\label{f}
\D f_t(\al,t)&=\frac{\rho^2-\rho^1}{2\pi}PV\int_{\R}
\frac{(\al-\beta)}{(\al-\beta)^2+(f(\al,t)-f(\beta,t))^2}(\partial_\al
f(\al,t)-\partial_\al f(\beta,t)) d\beta,
\\
f(\al,0)&=f_0(\al).
\end{split}
\end{align}
 The above equation can be linearized around the flat solution to find the following nonlocal partial differential equation
\begin{align*}
\begin{split}\label{le}
f_t(x,t)&=-\frac{\rho^2-\rho^1}{2}\Lambda f(x,t),\\
f(x,0)&=f_0(x),\quad x\in\R,
\end{split}
\end{align*}
where the operator $\Lambda$ is the square root of the Laplacian.
This linearization shows the parabolic character of the system.

Furthermore the stable system gives a maximum principle
$\|f\|_{L^{\infty}}(t)\leq \|f\|_{L^{\infty}}(0)$ \cite{DP2}; decay
rates are obtained for the periodic case:
 $$\|f\|_{L^\infty}(t)\leq \|f_0\|_{L^\infty}e^{-Ct},
 $$
and also for the case on the real line (flat at infinity):
$$
\D\|f\|_{L^\infty}(t)\leq \frac{\|f_0\|_{L^\infty}}{1+Ct}.
$$
There are several results on global existence for small initial
data (small compared to $1$  in several norms  more regular than
Lipschitz \cite{Peter, Yi2, SCH, DY, Esch2}) taking advantage of
the parabolic character of the equation for small initial data. In
\cite{ccgs} it is shown in the stable case that global existence for
solutions holds if the first derivative of the initial data is smaller
than an explicitly computable constant greater than $1/5$.
Furthermore, if $\|f_0\|_{L^\infty}<\infty$ and $\|\partial_x
f_0\|_{L^\infty}<1$, then there exists a global-in-time solution
that satisfies
$$
f(x,t)\in C([0,T]\times\R)\cap L^\infty([0,T];W^{1,\infty}(\R)),
$$
for each $T>0$. In particular $f$ is Lipschitz continuous.

Moreover, equation \eqref{f} yields an $L^2$ decay:
$$
\|f\|^2_{L^2}(t)+\frac{\rho^2-\rho^1}{2\pi}\!\int_0^t  ds ~\int_\R
d\al ~\int_\R dx~ \ln
\left(1+\Big(\frac{f(x,s)-f(\al,s)}{x-\al}\Big)^2\right)
\\
=\|f_0\|^2_{L^2},
$$
which does not imply, for large initial data, a gain of
derivatives in the system (see \cite{ccgs}). We will see below that the solutions to the Muskat problem with initial data in $H^4$ become real analytic immediately despite the weakness of the above decay formula.

 The main result we present here is:
\begin{thm}\label{oneone}
There exists a nonempty open set of initial data in $H^4$ with
Rayleigh-Taylor strictly positive $\sigma>0$ such that in finite
time the Rayleigh-Taylor $\sigma(\alpha,t)$ of the solution of \eqref{ec1d} is strictly
negative for all $\alpha$ in a nonempty open interval.
\end{thm}
 The geometry of this family of initial data is far from trivial:
numerical simulations performed in \cite{DPR} show that there exist
initial data with large steepness for which a regularizing  effect
appears. In fact, as  will be explained  in Section~\ref{csle}, the first
evidence of a change of sign in the Rayleigh-Taylor has
been experimentally found in a model with two interfaces.

 We proceed as follows:

 First, in section 3, we assume initial conditions at time $t=t_0$
that satisfy the Rayleigh-Taylor ($\sigma >0$) and the arc-chord
condition, and for which the boundary $z$ initially belongs to $H^4$. Let
$C_1$ be the constant in the arc-chord condition, let $C_2$ be an
upper bound for the $H^4$ norm of the initial data and let $c_3$ be a
lower bound for $\sigma$. Then there exists $t_1 > t_0$, with $t_1 $
depending only on $C_1, C_2, c_3$, such that the Muskat problem has a solution
for time $t\in [t_0,t_1]$, satisfying also the arc-chord and
Rayleigh-Taylor conditions. Moreover, for $t_0<t\leq t_1$, the
solution $ z(\alpha,t)$  is real analytic in a strip
$S(t)=\{\al+i\zeta: |\zeta|\leq c (t-t_0)\}$, where $c$ depends only on
$C_1, C_2, c_3$.

Our goal in section 4 is to show that the region of analyticity
does not collapse to the real axis as long as the Rayleigh-Taylor
is greater than or equal to 0. This allows us to reach a regime
for which the boundary $z$ develops a vertical tangent.

Section 5 is devoted to showing the existence of a large class of
analytic curves for which there exists a point where the tangent
vector is vertical and the velocities indicate that the curves are
going to turn over and reach the unstable regime for a small time.
Plugging these initial data into a Cauchy-Kowalewski theorem
indicates that the analytic curves turn over. Therefore the unstable
regime is reached.

 Finally, in section 6, a perturbative argument allows us to conclude that we can
find curves in $H^4$ close enough to the special class of analytic
curves described in Section 5, which satisfy the arc-chord and
Rayleigh-Taylor conditions. Then we can show the existence of the
curves passing the critical time and actually turning over. Therefore
the unstable regime is reached for an entire $H^4-$neighborhood of initial data.

\begin{rem}
In a forthcoming paper (see \cite{ADCPM3}) we will exhibit a particular initial datum for which we will show that once
the curve reaches the unstable regime the strip of analyticity
collapses in finite time and the solution breaks down. In section 8 we provide a very brief sketch of our proof of breakdown of smoothness
for the Muskat equation. These results were announced in \cite{turning}.
\end{rem}

\begin{rem}
The same approach can be done for the water waves problem, which
shows that, starting with some initial data given by
$(\al,f_0(\al))$, in finite time the interface reaches a regime in
which it is no longer a graph. Therefore there exists a time $t^*$
where the solution of the free boundary problem parametrized by
$(\al,f(\al,t))$ satisfies $\|f_\al\|_{L^\infty}(t^*)=\infty$ (see
section \ref{ww}). This scenario is known in the literature as
wave breaking \cite{CE} and there are numerical simulations
showing this phenomenon \cite{hou}.
\end{rem}

\begin{rem}
We conjecture that a result analogous to Theorem \ref{oneone} holds, in which surface tension is included.
We may simply use the same initial data as in Theorem \ref{oneone}, and take the coefficient of surface tension to be very small. The solutions are presumably changed only slightly by the surface tension (although we do not have a proof of this plausible assertion). Consequently, we believe that Muskat solutions with small surface tension can turn over.

A similar remark applies to water waves (see theorem
\ref{waterwaves}). There exist initial data for which water waves with surface tension turn over. A rigorous proof may be easily supplied, since local existence (backwards and forward in time) is known for water waves with surface tension
(see \cite{BL}).
\end{rem}

\section{The contour equation and numerical simulations}\label{csle}

Here we present the evolution equation in terms of the free boundary
which is going to be used throughout the paper, and the numerical
experiment that motivated the Theorem.

\subsection{The equation of motion}

By Darcy's law:

$$
\gradp\cdot u=-(\rho^2-\rho^1)\da z_2(\al)\delta(x-z(\al)),
$$
and Biot-Savart yields
\begin{equation}\label{ns}
z_t(\al)=-\frac{(\rho^2-\rho^1)}{2\pi}PV\int_{\R}\frac{(z(\al)-z(\al-\beta))^\bot}{|z(\al)-z(\al-\beta)|^2}\da
z_2(\al-\beta)d\beta.
\end{equation}
 For the first coordinate above one finds
$$
\frac{(\rho^2-\rho^1)}{2\pi}PV\int_{\R}\frac{(z_2(\al)-z_2(\al-\beta))}{|z(\al)-z(\al-\beta)|^2}\da
z_2(\al-\beta)d\beta
$$
$$=-\frac{(\rho^2-\rho^1)}{2\pi}PV\int_{\R}\frac{(z_1(\al)-z_1(\al-\beta))}{|z(\al)-z(\al-\beta)|^2}\da z_1(\al-\beta)d\beta
$$
using the identity
$$
PV\int_{\R}\partial_{\beta}\big(\ln (|z(\al)-z(\al-\beta)|^2)\big) d\beta=0.
$$
Therefore
$$
z_t(\al)=-\frac{(\rho^2-\rho^1)}{2\pi}PV\int_{\R}\frac{(z_1(\al)-z_1(\al-\beta))}{|z(\al)-z(\al-\beta)|^2}\da
z(\al-\beta)d\beta.
$$
Here we point out that in the Biot-Savart law the perpendicular
direction appears, but after the above integration by parts, we only see the tangential direction.

Adding the tangential term
$$
\frac{(\rho^2-\rho^1)}{2\pi}PV\int_{\R}\frac{(z_1(\al)-z_1(\al-\beta))}{|z(\al)-z(\al-\beta)|^2}d\beta \da z(\al),
$$
we find that the contour equation is given by
$$
z_t(\al)=\frac{(\rho^2-\rho^1)}{2\pi}PV\int_{\R}\frac{z_1(\al)-z_1(\al-\beta)}{|z(\al)-z(\al-\beta)|^2}(\da
z(\al)-\da z(\al-\beta))d\beta.
$$
For the $2\pi$ periodic interface the equation becomes
\begin{equation}\label{ecp}
z_t(\al)=\frac{(\rho^2-\rho^1)}{4\pi}\int_{-\pi}^{\pi}
\frac{\sin(z_1(\al)-z_1(\al-\beta))(\da z(\al)-\da
z(\al-\beta))}{\cosh(z_2(\al)-z_2(\al-\beta))-\cos(z_1(\al)-z_1(\al-\beta))}
d\beta.
\end{equation}
In order to see (\ref{ecp}) we take $z(\al)=z_1(\al)+iz_2(\al)$; it is easy to rewrite \eqref{ns}
as follows;
$$
\overline{z}_t(\al)=-\frac{(\rho^2-\rho^1)}{2\pi
i}PV\int_{\R}\frac{\da z_2(\beta)}{z(\al)-z(\beta)}d\beta.
$$
The classical identity
$$
\Big(\frac{1}{z}+\sum_{k\geq 1}\frac{z}{z^2-(2\pi
k)^2}\Big)=\frac{1}{2 \tan(z/2)}
$$
allows us to conclude that
$$
z_t(\al)=\frac{(\rho^2-\rho^1)}{4\pi}\int_{\T}
\frac{(\sinh(z_2(\al)-z_2(\beta)),-\sin(z_1(\al)-z_1(\beta)))}{\cosh(z_2(\al)-z_2(\beta))
-\cos(z_1(\al)-z_1(\beta))}\da z_2(\beta) d\beta,
$$
where $\T=\R/2\pi\Z$.

Analogously, using the equality
$$
\frac{(\rho^2-\rho^1)}{4\pi}PV\int_{\R}\frac{\sinh(z_2(\al)-z_2(\beta))}{\cosh(z_2(\al)-z_2(\beta))
-\cos(z_1(\al)-z_1(\beta))}\da z_2(\beta)d\beta
$$
$$
=-\frac{(\rho^2-\rho^1)}{4\pi}PV\int_{\R}\frac{\sin(z_1(\al)-z_1(\beta))}{\cosh(z_2(\al)-z_2(\beta))
-\cos(z_1(\al)-z_1(\beta))}\da z_1(\beta)d\beta
$$
and adding the appropriate tangential term, we obtain equation \eqref{ecp}.

\subsection{The scenario motivated by the numerics}

Our investigations started with the idea that interesting new phenomena may arise if we study three fluids, separated
by two interfaces. Careful numerical studies indicated that one of the interfaces may turn over. In attempting to prove
analytically the turnover indicated by the numerics, we discovered that a turnover can occur also for a single interface,
i.e., for the Muskat problem. This section describes one of our numerical experiments.

\begin{figure}
\centering
\includegraphics[width=.48\textwidth]{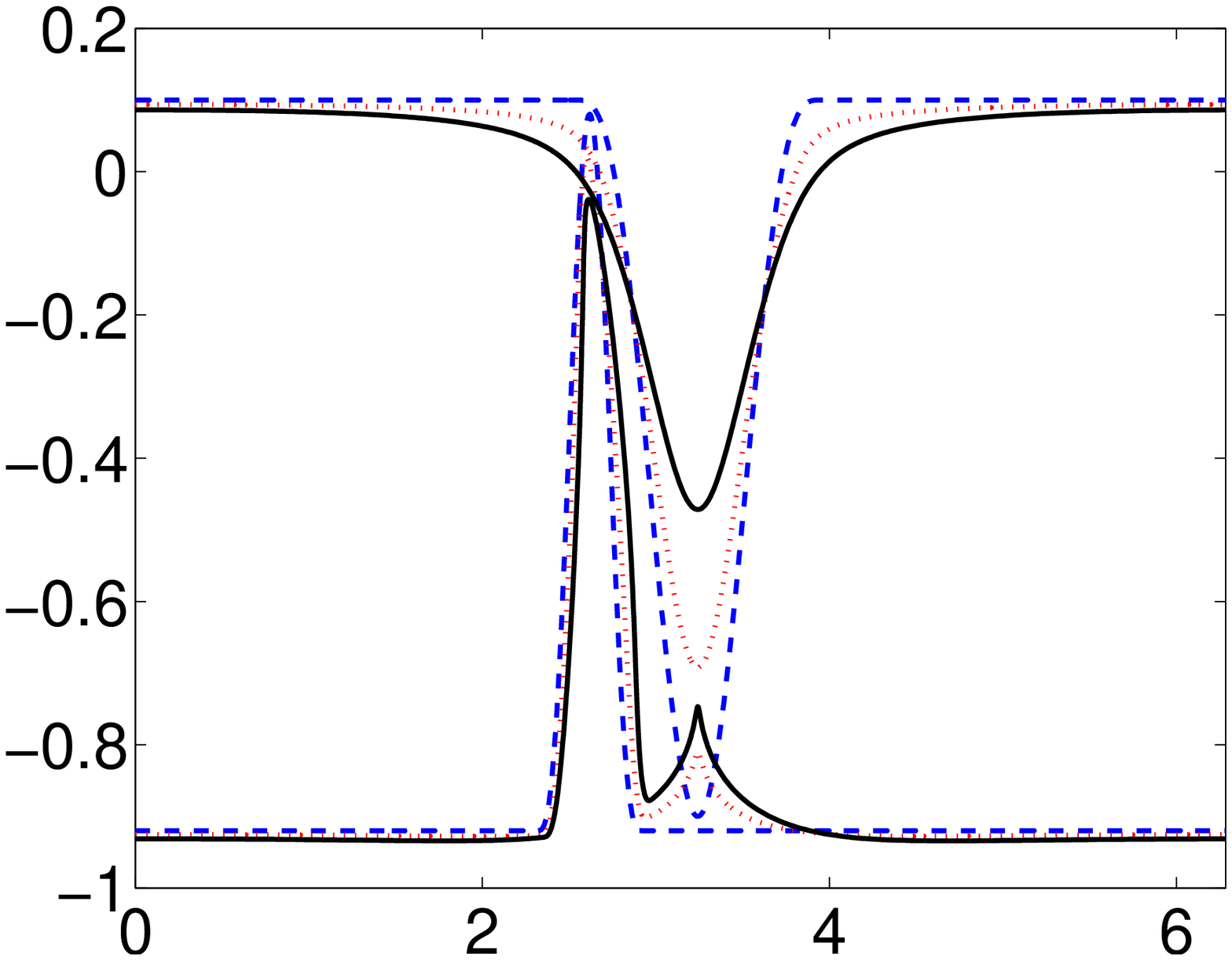}
\includegraphics[width=.48\textwidth]{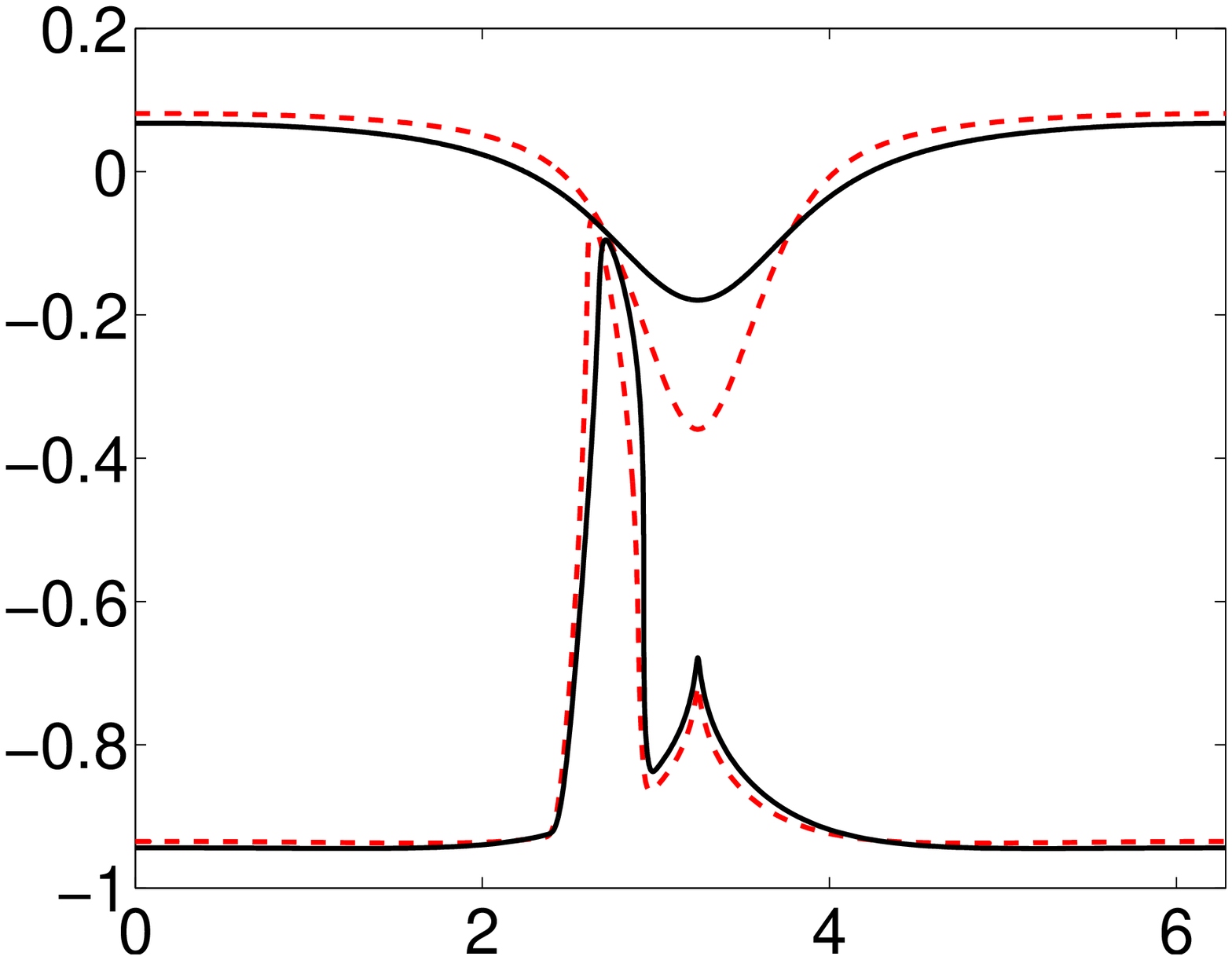}
\caption{{\em Left:} Solutions  to \eqref{eq2in} with initial data \eqref{inif}-\eqref{inig}
 at times $t=0$ (dashed blue), $t=3.46\cdot10^{-4}$ (red points) and $t=7.66\cdot 10^{-4}$ (black).
{\em Right:} Solutions at $t=1.04\cdot10^{-3}$ (dashed red) and $t=1.84\cdot10^{-3}$ (black)}\label{fig:sol}
\end{figure}

\begin{figure}
\centering
\includegraphics[width=.48\textwidth]{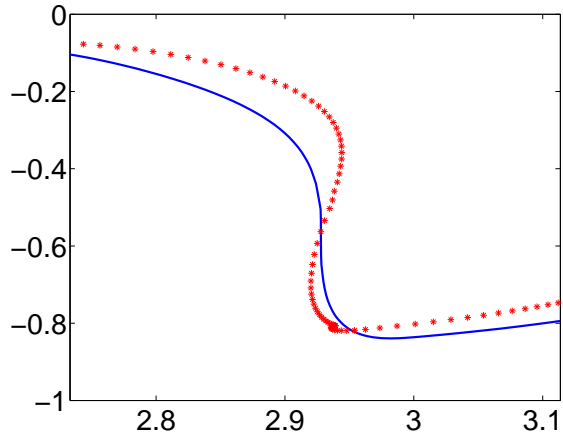}
\caption{ Zoom of the interface, modified so that its tangent is vertical at a  single point P;
 and the normal velocity along the curve, minus that at P, scaled by a  factor of 100.}\label{fig:campo}
\end{figure}

Proceeding as in the preceding section, one can derive the equations modeling the evolution of two interfaces separating three fluids with different densities $\rho_j$ ($j=1, 2,3$).
More precisely, assume that both
interfaces can be parametrized by  graphs $(\alpha,f(\alpha,t))$
and $(\alpha,g(\alpha,t))$, with $f$ lying above $g$. These
equations read in the periodic case, cf. \cite{DPR,DP3} (this
scenario has been recently also considered in \cite{Esch3}),
\begin{equation}\label{eq2in}
\begin{array}{ll}
f_t(\alpha,t) = \bar{\rho}_1 \, \mathcal{I}[f(\cdot,t),f(\cdot,t)] + \bar{\rho}_2\, \mathcal{I}[f(\cdot,t),g(\cdot,t)], \quad & f(\alpha,0)=f_0(\alpha),\\[.5em]
g_t(\alpha,t) = \bar{\rho}_2 \, \mathcal{I}[g(\cdot,t),g(\cdot,t)] + \bar{\rho}_1\, \mathcal{I}[g(\cdot,t),f(\cdot,t)], \quad & g(\alpha,0)=g_0(\alpha),
\end{array}
\end{equation}
where $\bar{\rho}_j=(\rho_{j+1}-\rho_j)/(4\pi)$, $j=1,2$, and, for given
functions $u(\alpha)$, $v(\alpha)$,
\begin{equation}\label{inteqs}
\mathcal{I}[u,v]:= PV \int_{\T} \frac{
(\da u(\alpha)-\da v(\alpha-\beta))
\tan(\beta/2)(1-\tanh^2( (u(\alpha)-v(\alpha-\beta))/2) )}
{\tan^2(\beta/2) + \tanh^2((u(\alpha)-v(\alpha-\beta))/2) } \, d\beta.
\end{equation}
The first terms $\mathcal{I}[f(\cdot,t),f(\cdot,t)]$ and
$\mathcal{I}[g(\cdot,t),g(\cdot,t)]$ in \eqref{eq2in} give the velocity of a
unique interface.
The cross terms $\mathcal{I}[f(\cdot,t),g(\cdot,t)]$ and
$\mathcal{I}[g(\cdot,t),f(\cdot,t)]$ take into account the interaction
of the two interfaces, and their contribution is getting bigger when the curves are getting closer.
This, together with the diffusive
behavior reported in \cite{DPR} for the equation
\begin{equation}\label{eq1in}
f_t(\alpha,t) = \bar{\rho}_1 \mathcal{I}[f(\cdot,t),f(\cdot,t)],\quad f(\alpha,0)=f_0(\alpha),
\end{equation}
and the mean conservation for $f$ and $g$, motivate the choice of
the following initial data, in the hope that some non regularizing
effect arises from the interaction of the two interfaces;
\begin{equation}\label{inif}
f_0(\alpha) = \left\{ \begin{array}{ll}
\displaystyle 0.1-\sin^3 \left(\frac{\pi(\alpha-M_1+r_1)}{2r_1} \right),
\quad & \mbox{ if }\ \alpha \in [M_1-r_1,M_1+r_1],\\[1em]
0.1, \quad & \mbox{otherwise}
\end{array}\right.
\end{equation}
and
\begin{equation}\label{inig}
g_0(\alpha) = \left\{ \begin{array}{ll}
\displaystyle  \sin^3 \left(\frac{\pi(\alpha-M_2+r_2)}{2r_2}
\right)^3-0.92, \quad & \mbox{ if }\ \alpha \in [M_2-r_2,M_2+r_2], \\[1em]
g_0(\alpha)=-0.92,\quad & \mbox{otherwise}.
\end{array} \right.
\end{equation}
The choice of  parameters $M_1=\pi+0.1$, $r_1 = 0.7$, $M_2 =
\pi/1.2$, $r_2 = 0.3$, $\bar{\rho}_1= 20\pi$ and $\bar{\rho}_2=
\pi/20$, yielded a strong growth of the derivative in the the lower
interface as the two curves approach, as shown in
Figure~\ref{fig:sol}.
Moreover, after introducing a small modification in the lower interface so that
the tangent at a certain point becomes actually infinite, and evaluating the
normal velocity relative to this point along the modified curve, we obtain the
result plotted in Figure~\ref{fig:campo}. This graphic clearly indicates that
the velocity field is forcing the interface to turn over.

The numerical approximation of \eqref{eq2in} addresses as a main
difficulty the absolute lack of knowledge about the behavior of
the solutions to \eqref{eq2in}. Indeed,
the goal of our experiments is precisely the search for some
singular behavior. The nonlocal terms make the computations
expensive and special care has to be taken in order to evaluate
the integrands in a neighborhood of $\beta=0$. For this, we used
Taylor expansions locally and computed exactly the principal
value. In this situation, adaptivity is strongly indicated, both
in space and time, since a good indicator of a singular behavior
will be given either by a sudden accumulation of spatial nodes or
a sudden reduction of the time steps.

In order to attain the highest resolution in the integration of
\eqref{eq2in} and compute the solutions shown in
Figure~\ref{fig:sol}, cubic spline interpolation of the curves
$f(\cdot,t)$ and $g(\cdot,t)$ with periodic boundary conditions was
used. This provides a $\mathcal{C}^2$ interpolant of each interface
at every time and allows, in particular, the evaluation of the
convolution terms at any $\beta \in [0,2\pi]$. Then, adaptive
quadrature can be applied to approximate the integrals and evaluate
the derivative at any time. In the experiments reported, adaptive
Lobatto quadrature was used, by means of the MATLAB routine {\tt
quadl}. For the time integration, the embedded Runge--Kutta formula
due to Dormand and Prince, DOPRI5(4), was implemented, since the
problem was not found to be particularly stiff, see for instance
\cite{Hai}. The time stepping was combined with a spatial node
redistribution after every successful step. For the redistribution
of the spatial nodes an algorithm following \cite{Drit} was
implemented, with some modifications taking into account that both
interfaces are graphs. For several tolerance requirements and
different choices of the parameters involved in the full adaptive
routine, the integration always failed  at a certain
critical time, suggesting the explosion of the derivative at a
certain point of the lower interface and the lack of validity of
\eqref{eq2in}, once this curve stops being a graph.

The phenomenon described above and the explicit representations of
the maximum of the solutions derived in \cite{DP2}, motivated the
search for special initial data which allowed us to understand
that this behavior also arises in the one-interface case.

\section{Instant Analyticity}

Here we show the main estimates that provide local-existence and
instant analyticity for a single curve that satisfies initially the
arc-chord and Rayleigh-Taylor conditions. We consider the function
$$
F(z)(\al,\beta)=\frac{\beta^2}{|z(\al)-z(\al-\beta)|^2},\quad
\al,\beta\in \R,
$$
and in the periodic setting
$$
F(z)(\al,\beta)=\frac{||\beta||^2}{2(\cosh(z_2(\al)-z_2(\al-\beta))
-\cos(z_1(\al)-z_1(\al-\beta)))},\quad \al,\beta\in \T,
$$
where $||x||= dist(x, 2\pi\Z)$.

If $F(z)\in L^\infty$ then we say that the curve satisfies the arc-chord condition, and the $L^{\infty}$ norm of $F$ is
called the arc-chord constant.

Let us clarify the meaning of the above arc-chord condition. Fix $t$, and assume that $z(\al,t)$ is a smooth
function of $\al$. Suppose $F\in L^\infty$. Letting $\beta$ tend to zero,
we conclude  that $|\da z(\al,t)|$  is bounded below. Since also $z$ is smooth, $|\da z(\al,t)|$ is also bounded above.
Consequently, the numerator in the fraction defining $F$
is comparable to the square of the arc-length between $z(\al,t)$ and $z(\al-\beta,t)$. On the other hand, the denominator of that
fraction is comparable to the square of the length of the chord joining $z(\al,t)$ to $z(\al-\beta,t)$. Thus, the boundedness of
$F$ expresses the standard arc-chord condition for the curve $z(\cdot,t)$ together with a lower bound for $|\da z(\al,t)|$.

\begin{thm}\label{instantanalycity}
Let $z(\al,0)=z_0(\al)\in H^4$, $F(z_0)(\al,\beta)\in L^\infty$ and
$\da z_1(\al,0)>0$ (R-T). Then there is a solution of the Muskat
problem $z(\al,t)$ defined for $0<t\leq T$ that continues
analytically into the strip $S(t)=\{\al+i\zeta: |\zeta|<ct\}$ for
each $t$. Here, $c$ and $T$ are determined by upper bounds of the
$H^4$ norm and the arc-chord constant of the initial data and a positive lower bound of $\partial_{\alpha}z_1(\alpha,0)$.
Moreover, for $0<t\leq T$, the quantity
$$
\sum_{\pm}\int(|z(\al\pm ict)-(\al+ict,0)|^2+|\da^4z(\al\pm ict)|^2) d\al
$$
is bounded by a constant determined by  upper bounds for the $H^4$
norm and the arc-chord constant of the initial data and a positive
lower bound of $\partial_{\alpha}z_1(\alpha,0)$. Above $|\cdot|$
is the modulus of a complex number or a vector in $\C^2$.
\end{thm}
Proof: For the proof we consider the contour $z\in H^4$ with
$z-(\al,0)$ periodic and $\da z_1(\al,0)>0$. In the case of the
real line similar arguments hold.  The Muskat equation reads
\begin{equation}\label{eqcp}
z_t(\al)=\int_{-\pi}^{\pi}
\frac{\sin(z_1(\al)-z_1(\al-\beta))(\da z(\al)-\da
z(\al-\beta))}{\cosh(z_2(\al)-z_2(\al-\beta))-\cos(z_1(\al)-z_1(\al-\beta))}
d\beta,
\end{equation}
where we suppose $\da z_1(\al,0)>0$. We also take
$\rho^2-\rho^1=4\pi$ since we are studying  the case
$\rho^2>\rho^1$. For the complex extension one finds
\begin{equation}\label{eqcomplex}
z_t(\al+i\zeta)=\int_{-\pi}^{\pi}
\frac{\sin(z_1(\al+i\zeta)-z_1(\al+i\zeta-\beta))(\da z(\al+i\zeta)-\da
z(\al+i\zeta-\beta))}{\cosh(z_2(\al+i\zeta)-z_2(\al+i\zeta-\beta))-\cos(z_1(\al+i\zeta)-z_1(\al+i\zeta-\beta))}
d\beta.
\end{equation}

We will use energy estimates. Consider
$$S(t)=\{\al+i\zeta\in\C:\,\al\in\T,\, |\zeta|<ct\},$$ for
$c$ given below\footnote{At the end of the proof we can take any
$c<\min_{\al}(\da z_1(\al,0)/|\da z(\al,0)|^2)$.},
$$
\|z\|^2_{L^2(S)}(t)=\sum_{\pm}\int_{\T}|z(\al\pm ict,t)-(\alpha\pm ict,0)|^2
d\al,
$$
$$
\|z\|^2_{H^k(S)}(t)=\|z\|^2_{L^2(S)}(t)+\sum_{\pm}\int_{\T}|\da^kz(\al\pm ict,t)|^2 d\al,
$$
where $k\geq 2$ as an integer, and
\begin{equation}\label{arc-chord-A}
F(z)(\al+i\zeta,\beta)=\frac{||\beta||^2}{2(\cosh(z_2(\al+i\zeta)-z_2(\al+i\zeta-\beta))
-\cos(z_1(\al+i\zeta)-z_1(\al+i\zeta-\beta)))},
\end{equation}
with norm
$$
\|F(z)\|_{L^\infty(S)}(t)=\sup_{\al+i\zeta\in
S(t),\beta\in\T}|F(z)(\al+i\zeta,\beta)|.
$$
Next, we define as follows:
$$
\|z\|^2_{S}(t)=\|z\|^2_{H^4(S)}(t)+\|F(z)\|_{L^\infty(S)}(t).
$$
We shall analyze the evolution of $\|z\|_{H^4(S)}(t)$.

 Before starting the energy estimates, we mention
an idea used previously e.g. in the proof of (6.3) in \cite{ADY}. Suppose $A(\al,\beta)$ is a $C^1(\T)$ function,
and suppose $f(\al)$ belongs to $L^2(\T)$. To estimate
\begin{equation}\label{discussion}
\int_{-\pi}^\pi A(\al,\al-\beta)\frac12\cot(\frac{\beta}{2})f(\al-\beta)d\beta
\end{equation}
we break up this integral as the sum of
\begin{equation}\label{fsh}
A(\al,\al)\int_{-\pi}^\pi \frac12\cot(\frac{\beta}{2})f(\al-\beta)d\beta
\end{equation}
and
\begin{equation}\label{fsc}
\int_{-\pi}^\pi \Big\{[A(\al,\al-\beta)-A(\al,\al)]\frac12\cot(\frac{\beta}{2})\Big\}f(\al-\beta)d\beta.
\end{equation}
The integral in \eqref{fsh} is simply the Hilbert transform of $f$ and the quantity in curly brackets in
\eqref{fsh} is bounded.
This idea will be used repeatedly, with $A(\al,\beta)$ arising from derivatives $\da^k z(\al,t)$ up to order 2,
and with $f(\al)=\da^4z_\mu(\al,t)$ ($\mu=1,2$).
Whenever we use this scheme, we will simply say that ``a Hilbert transform arises''. For similar simple ideas
used below, we refer the reader to the term $J_1$ in pg. 485 in \cite{ADY}.

Then, using above scheme, for the low order
terms in derivatives, it is easy to find that
\begin{equation}\label{npi}
\frac12\frac{d}{dt}\int_{\T}|z(\al\pm ict),t)-(\alpha\pm ic t,0)|^2d\al\leq C
(\|z\|_{S}(t)+1)^k.
\end{equation}
 In (\ref{npi}) and in several of the estimates, $k$
denotes a enough large  universal constant.

Next, we check that
$$
\frac12\frac{d}{dt}\int_{\T}|\da^4z(\al\pm ict,t)|^2d\al=
 \sum_{j=1,2}\frac12\frac{d}{dt}\int_{\T}|\da^4z_j(\al\pm ict,t)|^2d\al
$$
where
\begin{equation}\label{npiI}
\frac12\frac{d}{dt}\int_{\T}|\da^4z_j(\al\pm
ict,t)|^2d\al=\Re\int_\T \overline{\da^4z_j(\al\pm ict,t)}
(\partial_t(\da^4z_j)(\al\pm ict,t)\pm ic\da^5 z_j(\al\pm
ict,t))d\al.
\end{equation}

In order to simplify the exposition we write
$z(\al,t)=z(\al)$ for a fixed $t$, we treat both coordinates at the same time, we write $(x_1,x_2)\cdot(x_3,x_4)=x_1x_3+x_2x_4$ for $x_j\in\C$, $j=1,...,4$, we denote $\al\pm ict=\g$,
and we define
$$
Q(\g,\beta)=\cosh(z_2(\g)-z_2(\g-\beta))-\cos(z_1(\g)-z_1(\g-\beta)).
$$
Then we split the right hand side of (\ref{npiI}) by writing
$$
I_1=\Re\int_\T \overline{\da^4 z(\g)}\cdot \da^4z_t(\g)d\al,
$$
and
$$
I_2=\Re\int_\T \overline{\da^4 z(\g)}\cdot ic\da^5 z(\g)d\al.
$$
In $I_1$ we will find the R-T and use it to absorb $I_2$. We will
 decompose $I_1$ in order to find  the terms of at least fourth order. In order
to estimate the lower order terms, we refer the reader to the
paper \cite{ADY} (see, e.g., Lemma 6.1). We have $I_1=J_1+J_2+J_3+$ l.o.t., where
$$\|\text{l.o.t}\|_{L^2(\T)}\leq C(\|z\|_S+1)^k,$$
and $J_1,J_2,J_3$ are defined as follows:
\begin{align*}
J_1=\Re\int_\T \overline{\da^4 z(\g)}\cdot \Big(\int_\T
A(\g,\beta)
\frac{\da^4z_1(\g)-\da^4z_1(\g-\beta)}{Q(\g,\beta)}(\da
z(\g)-\da z(\g-\beta))d\beta \Big) d\al,
\end{align*}
where $A(\g,\beta)=\cos(z_1(\g)-z_1(\g-\beta))$,
$$
J_2=-\Re\int_\T \overline{\da^4 z(\g)}\cdot \Big(\int_\T \frac{
\sin(z_1(\g)- z_1(\g-\beta))}{(Q(\g,\beta))^2} (\da z(\g)-\da
z(\g-\beta))B(\g,\beta)d\beta \Big) d\al
$$
where
$$ B(\g,\beta)=(\sin(z_1(\g)-z_1(\g-\beta)),\sinh(z_2(\g)-z_2(\g-\beta)))
\cdot (\da^4 z(\g)-\da^4 z(\g-\beta)),$$
and
$$
J_3=\Re\int_\T \overline{\da^4 z(\g)}\cdot \Big(\int_\T \frac{\sin(
z_1(\g)- z_1(\g-\beta))}{Q(\g,\beta)}(\da^5 z(\g)-\da^5
z(\g-\beta))d\beta \Big) d\al.
$$
We split further $J_1=K_1+K_2$ where
$$
K_1=\Re\int_\T \overline{\da^4 z(\g)}\cdot
\da^4z_1(\g)\Big(PV\int_\T \frac{A(\g,\beta)}{Q(\g,\beta)}
(\da z(\g)-\da z(\g-\beta))d\beta\Big) d\al,
$$
$$
K_2=-\Re\int_\T \overline{\da^4 z(\g)}\cdot \Big(PV\int_\T \frac{A(\g,\beta)}{Q(\g,\beta)}(\da
z(\g)-\da z(\g-\beta))
\da^4z_1(\g-\beta)d\beta \Big) d\al.
$$
Taking into account the complex extension of the arc-chord condition, it is
easy to deal with $K_1$ to obtain
$$
K_1\leq (\|z\|_{S}(t)+1)^k.
$$
In $K_2$ it is possible to find a ``Hilbert transform'' applied to
$\da^4z_1$ as in (\ref{discussion}), and therefore an analogous estimate follows. We are done
with $J_1$. For $J_2$ we obtain similarly
$$
J_2\leq (\|z\|_{S}(t)+1)^k.
$$
Next, we split $J_3=K_3+K_4$ where
$$
K_3=\Re\int_\T \overline{\da^4 z(\g)}\cdot \da^5 z(\g)
\Big(PV\int_\T \frac{\sin(
z_1(\g)- z_1(\g-\beta))}{Q(\g,\beta)}d\beta \Big) d\al,
$$
$$
K_4=-\Re\int_\T \overline{\da^4 z(\g)}\cdot  \Big(PV\int_\T \frac{\sin(
z_1(\g)- z_1(\g-\beta))}{Q(\g,\beta)} \da^5
z(\g-\beta)d\beta \Big) d\al.
$$
We have to be careful, because $K_3$ for real curves is harmless, but for complex curves
we need to use the dissipative term to cancel out a dangerous term. We denote
\begin{equation}\label{df}
f(\g)=PV\int_\T \frac{\sin(
z_1(\g)- z_1(\g-\beta))}{Q(\g,\beta)}d\beta
\end{equation}
and therefore $K_3=L_1+L_2$ where
$$
L_1=\int_{\T}\Re(f)(\Re(\da^4z)\Re(\da^5z)+\Im(\da^4z)\Im(\da^5z))d\al,
$$
$$
L_2=\int_{\T}\Im(f)(-\Re(\da^4z)\Im(\da^5z)+\Im(\da^4z)\Re(\da^5z))d\al.
$$
An easy integration by parts allows us to get
$$
L_1=-\frac12\int_\T\Re(\da f)|\da^4z|^2d\al\leq  C(\|z\|_{S}(t)+1)^k.
$$
For $L_2$ we find
$$
L_2=\int_{\T}\Im(\da f)\Re(\da^4z)\Im(\da^4z)d\al +2\int_{\T}\Im(
f)\Im(\da^4z)\Re(\da^5z))d\al.
$$
The first term on the right is easy to dominate by
$C(||z||_{S}+1)^k$. We denote the second one by $M_1$. We claim
that
\begin{equation}\label{malo}
M_1\leq C(\|z\|_{S}(t)+1)^k+K\|\Im(f)\|_{H^2(S)} \|\la^{1/2}\da^4
z\|^2_{L^2(S)},
\end{equation}
for $K>0$ universal constant. To see this, we rewrite
$$
M_1=-2\int_{\T}\Im( f)\Im(\da^4z)\Re(\la(H(\da^4z)))d\al
$$
which yields
$$
M_1=-2\int_{\T}\la^{1/2}(\Im(f)\Im(\da^4z))\Re(\la^{1/2}(H(\da^4z)))d\al
$$
and therefore
\begin{align*}
M_1\leq&2\|\la^{1/2}(\Im(f)\Im(\da^4z))\|_{L^2(S)}\|\la^{1/2}\da^4z\|_{L^2(S)}\\
\leq&C\|\Im(f)\|_{H^2(S)}(\|\da^4z\|_{L^2(S)}+\|\la^{1/2}(\da^4z)\|_{L^2(S)})\|\la^{1/2}\da^4z\|_{L^2(S)}\\
\leq& C(\|z\|_{S}(t)+1)^k+K\|\Im(f)\|_{H^2(S)} \|\la^{1/2}\da^4
z\|^2_{L^2(S)}.
\end{align*}
Finally we find that
\begin{equation}\label{malo2}
K_3\leq C(\|z\|_{S}(t)+1)^k+K\|\Im(f)\|_{H^2(S)} \|\la^{1/2}\da^4
z\|^2_{L^2(S)}.
\end{equation}
 We will use the thickness of the strip to control the unbounded term above.

 For $K_4$ we decompose further: $K_4=L_3+L_4+L_5+L_6$ where
$$
L_3=-\Re \int_{-\pi}^\pi \overline{\da^4z(\g)}\cdot\int_{-\pi}^\pi
\frac{\beta^2}{Q(\g,\beta)}\frac1{\beta}
\big(\frac{\sin(z_1(\g)-z_1(\g-\beta))}{\beta}-\da z_1(\g)\big)
 \da^5z(\g-\beta)d\beta d\al,
$$
$$
L_4=-\Re\int_{-\pi}^\pi \overline{\da^4z(\g)}\cdot\da z_1(\g)\int_{-\pi}^\pi
\big(\frac{\beta^2}{Q(\g,\beta)}-\frac{2}{|\da
z(\g)|^2}\big)\frac{1}{\beta}\da^5 z(\g-\beta)d\beta d\al,
$$
$$
L_5=-\Re\int_{-\pi}^\pi\overline{\da^4z(\g)}\cdot\frac{\da z_1(\g)}{|\da
z(\g)|^2}\int_{-\pi}^\pi (\frac{2}{\beta}-\frac{1}{\tan (\beta/2)})\da^5
z(\g-\beta)d\beta d\al,
$$
$$
L_6=-\Re\int_{-\pi}^\pi\overline{\da^4z(\g)}\cdot\frac{\da z_1(\g)}{|\da
z(\g)|^2}\la(\da^4z)(\g) d\al.
$$
Inside $L_3$, $L_4$ and $L_5$ we can integrate by parts and
therefore
$$
L_3+L_4+L_5\leq C(\|z\|_S(t)+1)^k.
$$
In $L_6$ we use the splitting $L_6=M_2+M_3$ where
$$
M_2=\int_{\T}\Im(\frac{\da z_1}{|\da
z|^2})(-\Re(\da^4z)\cdot\Im(\la(\da^4z))+\Im(\da^4z)\cdot\Re(\la(\da^4z)))d\al,
$$
$$
M_3=-\int_{\T}\Re(\frac{\da z_1}{|\da
z|^2})(\Re(\da^4z)\cdot\Re(\la(\da^4z))+\Im(\da^4z)\cdot\Im(\la(\da^4z)))d\al.
$$
In $M_2$ it is easy to find a commutator formula:
$$
M_2=\int_{\T}[-\la\big(\Im(\frac{\da z_1}{|\da
z|^2})\Re(\da^4z)\big)+\Im(\frac{\da z_1}{|\da
z|^2})\Re(\la(\da^4z))]\cdot\Im(\da^4z)d\al,
$$
and the appropriate estimate follows. We find that $M_2\leq
C(||z||_S+1)^k$. For $M_3$ we write $M_3=N_1+N_2$ where
$$
N_1=-\int_{\T}[\Re(\frac{\da z_1}{|\da
z|^2})-m(t)](\Re(\da^4z)\cdot\Re(\la(\da^4z))+\Im(\da^4z)\cdot\Im(\la(\da^4z)))d\al,
$$
$$
N_2=-m(t)\|\la^{1/2}(\da^4 z)\|_{L^2(S)}^2,
$$
where
$$m(t)=\min_{\g}\Re(\frac{\da z_1(\g)}{|\da z(\g)|^2}).$$
We use  the pointwise estimate \cite{AD}
\begin{equation}\label{AD}
2g\la(g)-\la(g^2)\geq 0.
\end{equation}
Therefore
$$
N_1\leq \frac12\|\la(\Re(\frac{\da z_1}{|\da
z|^2}))\|_{L^\infty(S)}\|\da^4z\|^2_{L^2(S)}\leq C (\|z\|_S(t)+1)^k
$$
as long as
$$\Re(\frac{\da z_1(\g)}{|\da z(\g)|^2})>0.$$ Remember that initially $\Re(\frac{\da z_1(\g)}{|\da z(\g)|^2})$
is greater than zero (R-T). We will prove that it is going to keep like that for a short time.
For $I_2$ we find as before
$$
I_2=c\int_\T(\Im(\da^4z)(\g)\cdot\Re(\da^5z)(\g)-\Re(\da^4z)(\g)\cdot\Im(\da^5z)(\g))d\al\leq
c\|\la^{1/2}(\da^4 z)\|_{L^2(S)}^2.
$$
Finally
\begin{align*}
\frac12\frac{d}{dt}\int_{\T}|\da^4z(\al\pm ict)|^2d\al&\leq C
(\|z\|_S(t)+1)^k +(c+K\|\Im(f)\|_{H^2(S)}(t)-m(t))\|\la^{1/2}(\da^4
z)\|_{L^2(S)}^2(t).
\end{align*}
Note that $\|\Im(f)\|_{H^2(S)}(0)=0$.
If $c-m(0)<0$, we will show that
$$
c+K\|\Im(f)\|_{H^2(S)}(t)-m(t)<0
$$
for short time. It yields
$$
\frac12\frac{d}{dt}\int_{\T}|\da^4z(\al\pm ict)|^2d\al\leq C
(\|z\|_S(t)+1)^k,
$$
as long as $c+K\|\Im(f)\|_{H^2(S)}(t)-m(t)<0$. Using Sobolev
estimates, we proceed as in section 8 in \cite{ADY} to show that
$$
\frac{d}{dt}\|F(z)\|_{L^\infty(S)}\leq C(\|z\|_S(t)+1)^k.
$$
 From the two inequalities above and \eqref{npi} it is easy to obtain a priori energy estimates that depend
 upon the negativity of $c+K\|\Im(f)\|_{H^2(S)}(t)-m(t)$. We
get bona fide energy estimates as follows. We denote
$$
\|z\|_{RT}^2(t)=\|z\|^2_S(t)+1/(m(t)-c-K\|\Im(f)\|_{H^2(S)}(t)).
$$
At this point, it is easy to find that
$$
-\frac{d}{dt}\|\Im(f)\|_{H^2(S)}(t)\leq C(\|z\|_S(t)+1)^k
$$
using \eqref{df}, and therefore (see section 9 in \cite{ADY} for more details)
$$
\frac{d}{dt}\|z\|_{RT}(t)\leq C(\|z\|_{RT}(t)+1)^k.
$$
It follows that
$$
\|z\|_{RT}(t)\leq \frac{\|z\|_{RT}(0)+1}{(1-C(\|z\|_{RT}(0)+1)^kt)^{1/k}}-1,
$$
providing the a priori estimate with $C$ and $k$ universal constants.

 We approximate the problem as follows
\begin{eqnarray*} z_t^\ep(\al,t)&=&\phi_\ep*\int \frac{
\sin(\phi_\ep*z^\ep_1(\alpha)-\phi_\ep*z^\ep_1(\be))(\da(
\phi_\ep *z^\ep)(\alpha)-\pa_\be(\phi_\ep *z^\ep)(\beta))}{\cosh(z_2^\ep(\al)-z_2^\ep(\be))-\cos(z_1^\ep(\al)-z_1^\ep(\beta))}d\be\\
z^\ep(\al,0)&=&\phi_\ep * z_0(\alpha),
\end{eqnarray*}
where $\phi_\ep(x)=\phi(\al/\ep)/\ep$, $\phi$ is the heat kernel and
$\ep>0$. Picard's theorem yields the existence of a solution
$z^\ep(\al,t)$ in $C\left([0,T^\ep);H^4\right)$ which is analytic in
the whole space for $z_0$ satisfying the arc-chord condition and
$\ep$ small enough. Using the same techniques we have developed
above we obtain a bound for $z^\ep(\alpha,t)$ in $H^4$ in the strip
$S(t)$ for a small enough $T$ which is independent of $\ep$. We need
arc-chord, R-T, $z_0\in H^4$ and $c-m(0)<0$. Then we can pass to the
limit.

\section{Getting all the way to breakdown of
Rayleigh-Taylor}\label{getting} This section is devoted to proving
the following theorem.

\begin{thm}\label{RTmoe0} Let $z(\al,0)=z^0(\al)$ be an analytic
curve in the strip
$$S=\{\alpha+i\zeta\in \C\,:\, |\zeta| < h(0)\},$$
with $h(0)>0$ and
satisfying:
\begin{itemize}
\item The arc-chord condition, $F(z^0)(\al+i\zeta,\beta)\in
L^\infty(S\times \R)$ \item The Rayleigh-Taylor condition, $\da
z_1^0(\al)> 0$ . \item The curve $z^0(\alpha)$ is real for real
$\alpha$. \item The functions $z_1^0(\alpha)-\alpha$ and
$z_2^0(\alpha)$ are periodic with period $2\pi$. \item The
functions $z_1^0(\alpha)-\alpha$ and $z_2^0(\alpha)$ belong to
$H^4(\pa S)$. \end{itemize} Then there exist a time $T$ and a
solution of the Muskat problem $z(\al,t)$ defined for $0<t\leq T$
that continues analytically into some complex strip for each fixed
$t\in [0,T]$. Here $T$ is either a small constant depending only
on $||z^0||_{S}$ or it is the first time a vertical tangent
appears, whichever occurs first.
\end{thm}
Thus our Muskat solution is analytic as long as $\pa_\alpha
z_1(\al,t)\geq 0$.

We will use the following:
\begin{lemma}
Let $\varphi(\al\pm i\ze)=\sum_{k=-N}^{N}A_k e^{ik\al\mp k\ze}$.
Then, for $\ze>0$, we have
\begin{equation}
\frac{\partial}{\partial \ze}\sum_{\pm}\int_{\T}|\varphi(\al\pm
i\ze)|^2d\al\geq \frac{1}{10}\sum_{\pm}\int_\T \la\varphi(\al\pm
i\ze)\overline{\varphi(\al\pm i\ze)}d\al-10\int_{\T}
\la\varphi(\al)\overline{\varphi}(\al)d\al,
\end{equation}
where $\la\varphi(\al\pm i\ze)=\sum_{k=-N}^{N}|k|A_k
e^{ik\al}e^{\mp k\ze}$.
\end{lemma}

Proof: First we shall compute the left hand side in the frequency
space:
$$ \sum_{\pm}\int_{\T}|\varphi(\al\pm i\ze)|^2 d\al= 4\pi\sum_{k=-N}^N|A_k|^2\cosh(2|k|\ze).$$
On the other hand we have that
$$\sum_{\pm}\int_{\T}\Lambda \varphi(\al \pm i\ze)\overline{\varphi(\al\pm i\ze)}d\al=4\pi\sum_{k=-N}^N|k||A_k|^2\cosh(2|k|\ze),$$
while
$$\int_\T\Lambda \varphi(\al)\overline{\varphi(\al)} d\al=2\pi\sum_{k=-N}^N|k||A_k|^2.$$

Differentiating in $\ze$ we obtain

$$\frac{\partial}{\partial \ze}\int_{\T}|\varphi(\al\pm i\ze)|^2d\al=8\pi\sum_{k=-N}^{N}|k||A_k|^2\sinh(2|k|\ze).$$

The lemma holds since $\sinh(\ze)\geq\cosh(\ze)-1$ for any
$\ze>0$.

\begin{cor}\label{corolario}
Let $\varphi(\al\pm i\ze,t)=\sum_{k=-N}^{N}A_k(t)e^{ik\al}e^{\mp
k\ze}$ and $h(t)>0$ be a decreasing function of $t$. Then
$$
\frac{\partial}{\partial t}\sum_{\pm}\int_{\T}|\varphi(\al\pm i
h(t))|^2d\al\leq \frac{h'(t)}{10}\sum_{\pm}\int_\T
\la\varphi(\al\pm ih(t))\overline{\varphi(\al\pm ih(t))}d\al
$$
$$\qquad\qquad\qquad-10 h'(t)\int_{\T}\la\varphi(\al)\overline{\varphi}(\al)d\al+2\Re\sum_{\pm}\int_{\T}\varphi_t(\al\pm i h(t))
\overline{\varphi(\al\pm i h(t))}d\al.$$
\end{cor}

This corollary allows us to prove Theorem \ref{RTmoe0}.\\

Proof (Theorem \ref{RTmoe0}): The norms $\|z\|_{H^k(S)}$ and
$\|z\|_{S}$ are defined as before using the new strip $S(t)$
defined by

$$S(t)=\{\alpha+i\zeta\in \C\,:\, |\zeta|< h(t)\},$$
where $h(t)$ is a positive decreasing function of $t$.

We use the Galerkin approximation of  equation (\ref{eqcp}), i.e.
\begin{equation*}
\pa_t z^{[N]}(\zeta,t)=\Pi_N[J[z^{[N]}]](\zeta,t),\end{equation*}
where $\zeta\in \overline{S}(t)$, $\Pi_N$ will be specified below,
and
$$J[z](\alpha, t)=\int_{-\pi}^{\pi} \frac{\sin(z_1(\al)-z_1(\beta))(\da z(\al)-\da
z(\beta))}{\cosh(z_2(\al)-z_2(\beta))-\cos(z_1(\al)-z_1(\beta))}
d\beta.$$ We impose the initial condition
\begin{equation*}
z^{[N]}(\alpha,0)=z^{[N]}(\alpha).\end{equation*} Here, for a
large enough positive integer $N$, we define $z^{[N]}(\alpha,0)$
from $z^0(\alpha)$ by using the projection
\begin{equation*}
\Pi_N\,:\, \sum_{-\infty}^\infty A_k e^{ik\alpha}\mapsto
\sum_{-N}^N A_k e^{ik\alpha}.\end{equation*} We define
$z^{[N]}(\alpha)$ by stipulating that
\begin{equation*}
z_1^{[N]}(\alpha)-\alpha=\Pi_N[z_1^0(\alpha)-\alpha]\end{equation*}
and
\begin{equation*}
z_2^{[N]}(\alpha)=\Pi_N[z_2^0(\alpha)].\end{equation*} For $N$
large enough, the functions $z^{[N]}(\alpha,0)$ satisfy the
arc-chord and Rayleigh-Taylor condition.

We shall consider the evolution of the most singular quantity
$$\sum_{\pm}\int_{\T}|\da^4 z^{[N]}(\al\pm ih_N(t),t)|^2d\alpha,$$
where $h_N(t)$ is a smooth positive decreasing function on $t$,
with $h_N(0)=h(0)$, which will be given below. Also we denote
$$S_N(t)=\{\alpha+i\zeta\in \C\,:\, |\zeta|< h_N(t)\}.$$

>From now on, we will drop the dependency  on $N$ from $z^{[N]}$
and $h_N(t)$ in our notation. We will return to the previous notation in the discussion below
at the end of the section. Taking the derivative with respect to   $t$ yields

$$\frac{d}{dt}\int\limits_{\alpha\in\T}\left|\pa_{\alpha}^{4}z_\mu(\alpha\pm
ih(t),t)\right|^2d\alpha$$
$$=2\Re \int\limits_{\alpha\in\T} \overline{\pa_\alpha^{4}z_\mu(\alpha\pm ih(t),t)}\left\{\pa_t \pa_\alpha^{4}z_\mu(\alpha\pm ih(t),t)+ih'(t)\pa_\alpha^{5}z_\mu(\alpha\pm ih(t),t)\right\}d\alpha$$
$$=2\Re \int\limits_{\alpha\in\T} \overline{\pa_\alpha^{4}z_\mu(\alpha\pm ih(t),t)}\left\{\pa_\alpha^{4}\Pi_N[J_\mu[z]](\alpha\pm ih(t),t)+ih'(t)\pa_\alpha^{5}z_\mu(\alpha\pm ih(t),t)\right\}d\alpha$$
$$=2\Re \int\limits_{\alpha\in\T} \overline{\pa_\alpha^{4}z_\mu(\alpha\pm ih(t),t)}\left\{\Pi_N[\pa_\alpha^{4}J_\mu[z]](\alpha\pm ih(t),t)+ih'(t)\pa_\alpha^{5}z_\mu(\alpha\pm ih(t),t)\right\}d\alpha$$
\begin{equation*} =2\Re \int\limits_{\alpha\in\T}
\overline{\pa_\alpha^{4}z_\mu(\alpha\pm
ih(t),t)}\left\{\pa_\alpha^{4}J_\mu[z](\alpha\pm ih(t),t)
+ih'(t)\pa_\alpha^{5}z_\mu(\alpha\pm
ih(t),t)\right\}d\alpha,\end{equation*} since
$\pa_\alpha^4z_\mu(\alpha\pm ih(t),t)$ is a trigonometric
polynomial in the range of $\Pi_N$. Here $\mu=1$, $2$.

Using the above corollary we have that
$$\frac{d}{dt}\sum_{\pm}\int\limits_{\alpha\in\T}\left|\pa_{\alpha}^{4}z_\mu(\alpha\pm
ih(t),t)\right|^2d\alpha\leq \frac{h'(t)}{10}\sum_{\pm}\int_\T
\la(\da^4z_\mu)(\al\pm ih(t))\cdot\overline{\da^4 z_\mu(\al\pm
ih(t))}d\al$$
$$-10 h'(t)\int_{\T}\la(\da^4 z_\mu)(\al)\cdot\overline{\da^4z_\mu}(\al)d\al+2\sum_{\pm}\Re\int_{\T}
\da^4 J_\mu[z](\alpha,t)(\al\pm i
h(t))\cdot\overline{\da^4z_\mu(\al \pm ih(t))}d\alpha.$$ We shall
study in detail  the most singular term in $\pa^4 J[z](\alpha,t)$,
i.e. \begin{eqnarray*}\da ^ 4 J[z](\alpha\pm
ih(t),t)&=&\int_{-\pi}^\pi\frac{\sin(z_1(\al\pm
ih(t),t)-z_1(\be,t))(\da^5z(\al\pm
ih(t),t)-\pa_\beta^5z(\beta,t))}{\cosh(z_2(\al\pm
ih(t),t)-z_2(\be,t))-\cos(z_1(\al\pm ih(t),t)-z_1(\be,t))}
d\beta\\
\text{+ l.o.t}&\equiv& X + \text{l.o.t.,}
\end{eqnarray*}
where $||l.o.t||_{L^2(\T)}\leq C(||z||_S(t)+1)^k$ (see \cite{ADY} and our previous discussion of
 \eqref{discussion}).
We split $X$ in to the following terms
\begin{align*}
X  = & \int_{-\pi}^\pi K(\alpha\pm i h(t),\beta)(\da^5z(\al\pm
ih(t),t)-\pa_\beta^5z(\beta,t))d\be\\
& +\sigma(\alpha\pm
ih(t),t)\int_{-\pi}^\pi\cot\left(\frac{\alpha\pm
ih(t)-\beta}{2}\right)(\da^5z(\al\pm
ih(t),t)-\pa_\beta^5z(\beta,t))d\be\\
\equiv & X_1 + X_2,
\end{align*}
where
\begin{align*}
K(\alpha,\beta)  = &
\frac{\sin(z_1(\alpha,t)-z_1(\beta,t))}{\cosh(z_2(\alpha,t)-z_2(\beta,t))-\cos(z_1(\alpha,t)-z_1(\beta,t))}\\
& -\frac{\pa_\alpha z_1(\alpha,t)}{(\pa_\alpha
z_2(\alpha,t))^2+(\pa_\alpha
z_1(\alpha,t))^2}\cot\left(\frac{\alpha-\beta}{2}\right)
\end{align*}
and
\begin{align*}
\tilde{\sigma}(\alpha,t) = & \frac{\pa_\alpha z_1(\alpha,t)}{(\pa_\alpha
z_1(\alpha,t))^2+(\pa_\alpha z_2(\alpha,t))^2}.
\end{align*}
Let us denote $$\Gamma_{\pm}(t)=\{\zeta\in \C\,:\, \zeta=\alpha\pm i
h(t),\, \alpha\in\T\}.$$ Since $K(\alpha,\beta)$ is a holomorphic
function in $\alpha$ and $\beta$, with $\alpha,\beta\in S(t)$, for
fixed $t$ we  have that
\begin{align*}
X_1 =& \int_\pi^\pi K (\alpha\pm ih(t),\be)\da^5z(\alpha\pm
 ih(t),t)d\beta\\
& -\int_\pi^\pi K (\alpha\pm i h(t),\be)\da^5z(\beta,t)d\beta\\
\equiv & X_{11}+X_{12},
\end{align*}
and integration by parts shows that the term $X_{12}$ satisfies $||X_{12}||_{L^2(\T)}\leq C(||z||_{S}+1)^k$. In addition, we can write $X_{11}$ as follows
\begin{align*} X_{11} = &  \int_{w\in
\Gamma_{\pm}(t)}K(\alpha\pm i h(t),w)\pa^5z(\alpha\pm
ih(t),t)dw\\
= & P.V.\int_{w\in\Gamma_{\pm}(t)}\frac{\sin(z_1(\alpha\pm
ih(t),t)-z_1(w,t))\pa^5z(\alpha\pm ih(t),t)} {\cosh(z_2(\alpha\pm
ih(t),t)-z_2(w,t))-\cos(z_1(\alpha\pm ih(t),t)-z_1(w,t))}dw\\
& -\pa^5z(\alpha\pm ih(t),t)\sigma(\alpha\pm i
h(t),t)P.V.\int_{w\in
\Gamma_{\pm}(t)}\cot\left(\frac{\alpha\pm ih(t)-w}{2}\right)dw\\
= & P.V.\int_{-\pi}^\pi\frac{\sin(z_1(\alpha\pm
ih(t),t)-z_1(\beta\pm ih(t),t))\pa^5z(\alpha\pm ih(t),t)}
{\cosh(z_2(\alpha\pm ih(t),t)-z_2(\beta\pm i
h(t),t))-\cos(z_1(\alpha\pm ih(t),t)-z_1(\beta\pm i h(t),t))}d\be,
\end{align*}
As before we call \begin{align*}&f(\alpha\pm i h(t),t)\\
&= P.V.\int_{-\pi}^\pi\frac{\sin(z_1(\alpha\pm
ih(t),t)-z_1(\beta\pm ih(t),t))} {\cosh(z_2(\alpha\pm
ih(t),t)-z_2(\beta\pm i h(t),t))-\cos(z_1(\alpha\pm
ih(t),t)-z_1(\beta\pm i h(t),t))}d\be\\
&=P.V.\int_{-\pi}^\pi\frac{\sin(z_1(\alpha\pm
ih(t),t)-z_1(\alpha\pm ih(t)-\beta,t))} {\cosh(z_2(\alpha\pm
ih(t),t)-z_2(\alpha\pm i h(t)-\beta,t))-\cos(z_1(\alpha\pm
ih(t),t)-z_1(\alpha\pm i h(t)-\beta,t))}d\be.
\end{align*} Thus
$$X_{11}=\pa^5z(\alpha\pm ih(t),t)f(\alpha\pm i h(t),t).$$
Also we can write $X_2$ in the following way;
\begin{align*}
X_2 = & \tilde{\sigma}(\alpha\pm
ih(t),t)\int_{-\pi}^\pi\cot\left(\frac{\alpha\pm
ih(t)-\beta}{2}\right)(\da^5z(\al\pm
ih(t),t)-\pa_\beta^5z(\beta,t))d\be\\
=& \tilde{\sigma}(\alpha\pm
ih(t),t)\int_{w\in\Gamma_\pm(t)}\cot\left(\frac{\alpha\pm
ih(t)-w}{2}\right)(\da^5z(\al\pm
ih(t),t)-\pa_\beta^5z(w,t))dw
\end{align*}
\begin{align*}
=& \tilde{\sigma}(\alpha\pm
ih(t),t)P.V.\int_{w\in\Gamma_\pm(t)}\cot\left(\frac{\alpha\pm
ih(t)-w}{2}\right)\da^5z(\al\pm
ih(t),t)dw\\
&- \tilde{\sigma}(\alpha\pm
ih(t),t)P.V.\int_{w\in\Gamma_\pm(t)}\cot\left(\frac{\alpha\pm
ih(t)-w}{2}\right)\da^5z(w,t)dw
\end{align*}
\begin{align*}
= & -\tilde{\sigma}(\alpha\pm
ih(t),t)P.V.\int_{-\pi}^\pi\cot\left(\frac{\alpha-\beta}{2}\right)\da^5z(\beta\pm
ih(t),t)d\beta\\
= & -\tilde{\sigma}(\alpha\pm
ih(t),t)P.V.\int_{-\pi}^\pi\frac{1}{2}\csc^2\left(\frac{\alpha-\beta}{2}\right)(\da^4z(\alpha\pm
i h(t),t)-\da^4z(\beta\pm ih(t),t))d\beta
\end{align*}
and finally
\begin{align*}
X_2= & -2\pi\tilde{\sigma}(\alpha\pm ih(t),t)(\Lambda \da^4 z) (\alpha\pm i
h(t),t).
\end{align*}
Then we find two dangerous terms
$$I_1=2\Re\int_{\T}f(\al\pm i h(t),t)\overline{(\da^4z_\mu)(\al\pm i
h(t))}\cdot(\da^5z_\mu)(\al\pm i h(t))d\al$$ and $$
I_2=-4\pi\Re\int_{\T}\tilde{\sigma}(\alpha\pm i h(t),t)\la (\da^4
z_\mu)(\al\pm i h(t))\cdot\overline{\da^4z_\mu(\al \pm
ih(t))}d\alpha.$$ The rest can be bounded by $C(\|z\|_{S}+1)^k(t)$
as in the previous section. In order to bound $I_1$ and $I_2$ we
use the following commutator estimate:

\begin{eqnarray}
||\Lambda^{\frac{1}{2}}(fg)-f\Lambda^{\frac{1}{2}}g||_{L^2(\T)}
\leq C||\Lambda^{1+\ep}f||_{L^2(\T)}||g||_{L^2(\T)},\label{commutator}
\end{eqnarray}
for $f(\alpha)=\sum_{-N}^{N}f_k e^{ikx}$ and
$g(\alpha)=\sum_{-N}^{N}g_k e^{ikx}$, where $\ep>0$ and $C$ does
not depend on $N$. The proof of (\ref{commutator}) will be left to
the reader.

First we estimate $I_1$. We denote $\gamma=\alpha+i h(t)$.
\begin{align*}
I_1 = & 2\Re \int_{-\pi}^\pi f(\gamma,t)\overline{\da^4
z_\mu(\gamma,t)}\da^5 z_\mu(\gamma,t) d\al\\
= & 2\int_{-\pi}^\pi
\Re(f(\gamma))\left\{\Re(\da^4z_\mu(\gamma,t))\pa_\alpha(\Re(\da^4z_\mu(\gamma,t)))+
\Im(\da^4z_\mu(\gamma,t))\pa_\alpha(\Im(\da^4z_\mu(\gamma,t)))\right\}d\alpha\\
- & 2\int_{-\pi}^\pi
\Im(f(\gamma))\left\{\Re(\da^4z_\mu(\gamma,t))\pa_\alpha(\Im(\da^4z_\mu(\gamma,t)))+
\Im(\da^4z_\mu(\gamma,t))\pa_\alpha(\Re(\da^4z_\mu(\gamma,t)))\right\}d\alpha\\
\equiv & I_{11}+I_{12}.
\end{align*}
Integrating by parts we have that $||I_{11}||_{L^2(\T)}\leq C(||z||_{S}+1)^k$. In
order to estimate $I_{12}$ we note that $f(\gamma,t)$ is real for
real $\gamma$. Then
$$\Im (f(\alpha\pm i h(t),t))=h(t)\tilde{f}_{\pm}(\alpha,t),$$
where $$||\tilde{f}_{\pm}||_{H^2(\T)}\leq C(||z||_S(t)+1)^k.$$
Then we can write
\begin{align*}
\int_{-\pi}^\pi
&\Im(f(\gamma))\Re(\da^4z_\mu(\gamma,t))\pa_\alpha(\Im(\da^4z_\mu(\gamma,t)))d\alpha\\
&=h(t)\int_{-\pi}^\pi
\tilde{f}_\pm(\alpha,t)\Re(\da^4z_\mu(\gamma,t))\pa_\alpha(\Im(\da^4z_\mu(\gamma,t)))d\alpha\\
&=-h(t)\int_{-\pi}^\pi
\tilde{f}_\pm(\alpha,t)\Re(\da^4z_\mu(\gamma,t))\Lambda H(\Im(\da^4z_\mu(\gamma,t)))d\alpha\\
&=-h(t)\int_{-\pi}^\pi
\Lambda^\frac{1}{2}(\tilde{f}_\pm(\alpha,t)\Re(\da^4z_\mu)(\gamma,t))\Lambda^{\frac{1}{2}}
 H(\Im(\da^4z_\mu(\gamma,t)))d\alpha\\
&=-h(t)\int_{-\pi}^\pi
\left\{\Lambda^\frac{1}{2}(\tilde{f}_\pm(\alpha,t)\Re(\da^4z_\mu)(\gamma,t))-
\tilde{f}_\pm(\alpha)\Lambda^{\frac{1}{2}}\Re(\da^4z_\mu)\right\}
\Lambda^{\frac{1}{2}} H(\Im(\da^4z_\mu(\gamma,t)))d\alpha\\
&-h(t)\int_{-\pi}^{\pi}\tilde{f}_\pm(\alpha)\Lambda^{\frac{1}{2}}\Re(\da^4z_\mu)\Lambda^{\frac{1}{2}}
H(\Im(\da^4z_\mu(\gamma,t)))d\alpha\\
& \leq
h(t)||\Lambda^\frac{1}{2}(\tilde{f}_\pm(\cdot,t))\Re(\da^4z_\mu(\cdot\pm
i h(t),t))-
\tilde{f}_\pm(\cdot,t)\Lambda^{\frac{1}{2}}\Re(\da^4z_\mu(\cdot\pm
i
h(t)))||_{L^2(\T)}\\
&\times||\Lambda^{\frac{1}{2}} H(\Im(\da^4z_\mu(\cdot\pm
i h(t),t)))||_{L^2(\T)}\\
&+h(t)||\tilde{f}_\pm||_{L^\infty(\T)}||\Lambda^\frac{1}{2}\Re(\da^4z_\mu(\cdot\pm
i h(t)))||_{L^2(\T)} ||\Lambda^\frac{1}{2}H\Im(\da^4z_\mu(\cdot\pm
i h(t)))||_{L^2(\T)}.
\end{align*}
Using  the estimate (\ref{commutator}) yields
\begin{align*}
&\int_{-\pi}^\pi
\Im(f(\gamma))\Re(\da^4z_\mu(\gamma,t))\pa_\alpha(\Im(\da^4z_\mu(\gamma,t)))d\alpha\\
&\leq
h(t)||\Lambda^{1+\ep}\tilde{f}_\pm||_{L^2(\T)}||\Re(\da^4z_\mu(\cdot\pm
i h(t),t))||_{L^2(\T)}||\Lambda^{\frac{1}{2}}
(\Im(\da^4z_\mu(\cdot\pm
i h(t),t)))||_{L^2(\T)}\\
&+h(t)||\tilde{f}_\pm||_{L^\infty(\T)}||\Lambda^\frac{1}{2}\Re(\da^4z_\mu(\cdot\pm
i h(t)))||_{L^2(\T)} ||\Lambda^\frac{1}{2}\Im(\da^4z_\mu(\cdot\pm
i h(t))||_{L^2(\T)}\\
&\leq C h(t)(||z||_S+1)^k+ Ch(t)(||z||_S+1)^k||\Lambda^\frac{1}{2}\da^4z_\mu(\cdot\pm
i h(t),t)||^2_{L^2(\T)}\\
&= Ch(t)(||z||_S+1)^k+ Ch(t)(||z||_S+1)^k\int_{-\pi}^\pi
\overline{\da^4z_\mu(\gamma,t)}\Lambda\da^4z_\mu(\gamma,t)d\alpha.
\end{align*}
Now $I_1$ is equal to the integral to the left, plus a similar
integral that can be bounded in a similar way.

Thus we obtain that
\begin{equation}\label{aaf}
\sum_{\pm} I_1\leq C(\|z\|_{S}+1)^k+Ch(t)(\|z\|_{S}+1)^k \|\la^{1/2}\da^4
z\|^2_{L^2(S)}.
\end{equation}

By assumption the R-T $\tilde{\sigma}$ is bigger than zero for
real values. In order to avoid problems with the imaginary part we
may write $$\frac{\da z_1(\al\pm i h(t),t)}{(\da z_1(\al\pm i
h(t)))^2+(\da z_2(\al\pm i h(t)))^2}=\frac{\da z_1(\al,t)}{|\da
z(\al,t)|^2}+h(t)g_{\pm}(\alpha,t).$$ where
$$||g_\pm||_{H^2(\T)} \leq C(||z||_{S}+1)^k.$$

One finds,

$$I_2=-2\Re\int_{\T}\frac{\da z_1(\al)}{|\da z(\al)|^2}\la (\da^4 z_{\mu})(\al\pm i h(t))\cdot\overline{\da^4z_\mu(\al \pm ih(t))}d\alpha$$
$$-h(t)2\Re\int_{\T}g_{\pm}(\al,t)\la (\da^4 z_\mu)(\al\pm i h(t))\cdot\overline{\da^4 z_\mu(\al\pm ih(t))}d\alpha.$$

The first term above can be treated as  in  section 3 taking
advantage of the inequality \eqref{AD}. Here we just need $\da
z_1(\al)\geq 0$. The second term can be treated using the
inequality (\ref{commutator}) as  with the term $I_1$. We find that
$$
\sum_{\pm} I_2\leq
C(||z||_{S}+1)^k+Ch(t)\|g_{\pm}\|_{H^2(S)}\|\la^{1/2}\da^4
z\|^2_{L^2(S)},
$$
and therefore
\begin{equation}\label{aag}
\sum_{\pm} I_2\leq C(||z||_{S}+1)^k+Ch(t)(\|z\|_{S}+1)^k \|\la^{1/2}\da^4
z\|^2_{L^2(S)}.
\end{equation}

Using \eqref{aaf} and \eqref{aag} we have that
$$
\frac{d}{dt}\sum_{\pm}\int_{\T}|\da^4 z_\mu(\al\pm ih(t))|^2
d\alpha\leq C(\|z\|_{S}(t)+1)^k-10 h'(t)\int_{\T}\la(\da^4
z_\mu)(\al)\cdot\overline{\da^4z_\mu}(\al)d\al$$
$$
+(C(\|z\|_{S}(t)+1)^k h(t)+ \frac{1}{10}h'(t))\int_\T
\la(\da^4z_\mu)(\al\pm ih(t))\cdot\overline{\da^4 z_\mu(\al\pm
ih(t))}d\al.
$$

Choosing
$$h(t)=h(0)\exp(-10C\int_0^t(\|z\|_{S}+1)^k(r)dr)$$
we eliminate the most dangerous term. The other term in the
expression above involves with a function on the real line and it is
easily controlled. Indeed
$$\int_{\T}\la \da^4 z_\mu(\al)\cdot\da^4z_\mu(\al)\leq\frac{C}{h(t)}\sum_{\pm}\int_{\T}|\da^4 z_\mu(\alpha\pm ih(t))|^2d\alpha,$$
as one sees by examining the Fourier expansion of $\da^4
z_\mu(\al,t)$.

Thus
$$\left|10 h'(t)\int_{\T}\la(\da^4
z_\mu)(\al)\cdot\da^4z_\mu(\al)d\al\right|\leq
C\frac{|h'(t)|}{h(t)}||z||^2_{S}\leq C(||z||_S+1)^{k+2}.$$

And we obtain finally
$$
\frac{d}{dt}\sum_{\pm}\int_{\T}|\da^4 z(\al\pm ih(t))|^2
d\alpha\leq C(\|z\|_{S}(t)+1)^{k+2}.
$$
Recovering the dependency on $N$ in our notation we have that
\begin{equation}\label{TN}
\frac{d}{dt}\sum_{\pm}\int_{\T}|\da^4 z^{[N]}(\al\pm ih_N(t))|^2
d\alpha\leq C(\|z^{[N]}\|_{S_N}(t)+1)^{k+2}.
\end{equation}

As  in the previous section, we can obtain a bound of the
evolution of the arc-chord condition that depends on
$C(\|z^{[N]}\|_{S_N}(t)+1)^{k+2}$.

This estimate is true whenever $t\in [0,T_N]$, where $T_N$ is the
maximal time of existence of the solution $z^{[N]}$. In addition
inequality (\ref{TN}) shows that we can extend these solutions in
$H^4(S)$ up to a small enough time $T$ independent of $N$ and
depending on the initial data.

The above calculation shows that the strip may shrink but does not
collapse as long as $\da z_1(\al,t)\geq 0$.

\section{From an analytic curve in the stable regime to an analytic curve in the unstable
regime}\label{AsAu} In this section we show that there exist some
initial data which are analytic curves satisfying the arc-chord
and R-T conditions such that the
solution of the Muskat
  problem reaches the unstable regime. In order to do it we will prove the local existence of solutions for analytic
  initial data without assuming the  R-T condition. Then we will construct some suitable initial data for our purpose.

\begin{thm}\label{kovalevsky}
Let $z_0$ be an analytic curve satisfying the arc-chord condition.
Then there exists an analytic solution for the Muskat problem
in some interval $[-T,T]$ for a small enough $T>0$.
\end{thm}
\begin{rem} Notice that in theorem (\ref{kovalevsky}) there is no assumption on the R-T condition. The proof we use here
is analogous to the one in \cite{SSBF} based on  Cauchy-Kowalewski theorems \cite{Nirenberg,Nishida}
(for an application to the Euler equation see \cite{Bardos}). Here we cannot parametrize the curve as a graph,
so we have to change the argument substantially in the proof in order to deal with the arc-chord condition.\end{rem}

Proof: We use the same notation as before. Let
$\{X_r\}_{r>0}$ be a scale of Banach spaces given by $\R^2-$valued
real functions $f$ that can be extended  into the complex
strip $S_r=\{\al+i\zeta\in\C: |\zeta|< r\}$ such that the  norm
$$
\|f\|^2_r=\sum_{\pm}\int_\T|f(\al\pm ir)-(\al\pm
ir,0)|^2d\al+\int_\T|\da^4f(\al\pm ir)|^2d\al,
$$
is finite and $f(\alpha)-(\alpha,0)$ is $2\pi-$periodic.

Let $z^0(\al)$ be a curve satisfying the arc-chord condition and
$z^0(\al)\in X_{r_0}$ for some $r_0>0$. Then, we will show that
there exist a time $T>0$ and $0<r<r_0$ so that there is a unique
solution to \eqref{eqcomplex} in $C([0,T];X_r)$.

It is easy to check that $X_{r}\subset X_{r'}$ for $r'\leq r$ due to
the fact that $\|f\|_{r'}\leq \|f\|_{r}$. A simple application of
the Cauchy formula gives
\begin{equation}\label{Cauchy}
\|\da f\|_{r'}\leq \frac{C}{r-r'}\|f\|_{r},
\end{equation}
for $r'<r$. Next, we write equation \eqref{eqcomplex} as follows:
$$z_t(\al+i\zeta,t)=G(z(\al+i\zeta,t)),$$ with
$$G(z(\al+i\zeta,t))=\int_{-\pi}^{\pi}
\frac{\sin(z_1(\al+i\zeta)-z_1(\al+i\zeta-\beta))(\da z(\al+i\zeta)-\da
z(\al+i\zeta-\beta))}{\cosh(z_2(\al+i\zeta)-z_2(\al+i\zeta-\beta))-\cos(z_1(\al+i\zeta)-z_1(\al+i\zeta-\beta))}
d\beta.$$

We take $0\leq r'<r$ and we introduce the open set $O$ in $S_{r}$ given by
\begin{equation}\label{openO}
O=\{z,\omega\in X_{r}: \|z\|_{r}<R,\quad
\|F(z)\|_{L^\infty(S_r)}<R^2\},
\end{equation}
with $F(z)(\al+i\zeta,\beta,t)$ given by \eqref{arc-chord-A}.
Then the function $G$ for $G:O\rightarrow X_{r'}$ is a continuous
mapping. In addition, there is a constant $C_R$ (depending on $R$
only) such that
\begin{equation}\label{cota}
\|G(z)\|_{r'}\leq \frac{C_R}{r-r'}\|z\|_{r},
\end{equation}
\begin{equation}\label{casiL}
 \|G(z^2)-G(z^1)\|_{r'}\leq \frac{C_R}{r-r'}\|z^2-z^1\|_{r},
\end{equation}
and
\begin{equation}\label{paraarc-chord}
\sup_{\al+i\zeta\in S_r,\beta\in\T}
|G(z)(\al+i\zeta)-G(z)(\al+i\zeta-\beta)|\leq C_R|\beta|,
\end{equation}
for $z,z^j\in O$. The above inequalities can be proved by estimating as in previous sections. Then they yield the proof of
theorem \ref{kovalevsky}. The argument is analogous to
\cite{Nirenberg} and \cite{Nishida}. We have to deal with the
arc-chord condition so we will point out the main differences. For
initial data $z^0\in X_{r_0}$ satisfying arc-chord, we can
find a $0<r_0'<r_0$ and a constant $R_0$ such that $\|z^0\|_{r_0'}<
R_0$ and
\begin{equation}\label{arc-chord-}
2\frac{\cosh(z^0_2(\al+i\zeta)-
z^0_2(\al+i\zeta-\beta))-\cos(z^0_1(\al+i\zeta)-z^0_1(\al+i\zeta-\beta))}{||\beta||^2}>\frac{1}{R_0^2},
\end{equation}
for $\al+i\zeta\in S_{r_0'}$. We take $0<r<r_0'$ and $R_0<R$ to
define the open set $O$ as in \eqref{openO}. Therefore we can use
the classical method of successive approximations:
\begin{equation}\label{iteracion}
z^{n+1}(t)=z^0+\int_0^t G(z^n(s))ds,
\end{equation}
for $G:O\rightarrow X_{r'}$ and $0< r'<r$. We assume by induction
that $$\|z^k\|_r(t)< R, \qquad\mbox{ and }\qquad
\|F(z^k)\|_{L^\infty(S_r)}(t)< R$$ for $k\leq n$ and $0<t<T$ with
$T=\min(T_A,T_{CK})$ and $T_{CK}$ the time obtained in the proofs in
\cite{Nirenberg} and \cite{Nishida}, and $T_A$ determined below.
Now, we will check that $\|F(z^{n+1})\|_{L^\infty(S_r)}(t)< R$
for suitable $T_A$. The rest of the proof follows in the same way as in
\cite{Nirenberg}, \cite{Nishida}.

Definitions (\ref{iteracion}) and (\ref{arc-chord-A}) easily imply that
\begin{align*}
\begin{split}
|(F(z^{n+1})(\al+i\zeta,\beta,t))^{-1}|&\geq
|(F(z^0)(\al+i\zeta,\beta,t))^{-1}|-C_R(t^2+t)\geq
\frac{1}{R_0^2}-C_R(t^2+t).
\end{split}
\end{align*}
To see this, we just use the formulas for $\cos(a+b)$ and
$\cosh(a+b)$, and bounds for the functions
$\frac{\cosh(x)-1}{x^2}$, $\frac{1-\cos(x)}{x^2}$,
$\frac{\sinh(x)}{x}$, $\frac{\sin(x)}{x}$, for bounded $x$.
Therefore, taking $$0<T_A< \min
\left\{1,\sqrt{\left(\frac{1}{R_0^2}-\frac{1}{R^2}\right)\frac{1}{2C_{R}}}\right\},$$
we obtain $\|F(z^{n+1})\|_{L^\infty(S_r)}(t)< R$. This completes the proof of theorem \ref{kovalevsky}.\\

The next step will be the construction of analytic initial data
such that
$$
\begin{array}{ll}
a.\,\,\da z_1(\alpha)>0\mbox{ if }\alpha \neq 0.\qquad & b.\,\,\da z_1(0)=0.\\
&\\
c.\,\,\da z_2(0)>0. & d.\,\, \da v_1(0)<0.\\
\end{array}
$$
Also $z_1(\alpha)-\alpha$ and $z_2(\alpha)$ are
$2\pi-$periodic.

Here $v_\mu(\alpha,t)$, with $\mu=1,2$, are the velocities given by
$$v_\mu(\alpha,t)=\int_{-\pi}^\pi\frac{\sin(z_1(\al)-z_1(\be))}{\cosh(z_2(\alpha)-z_2(\beta))-\cos(z_1(\alpha)-z_1(\beta))}(\da
z_\mu(\al)-\da z_\mu(\be))d\be.$$

Notice that in this situation the graph $f : \ff{R}\rightarrow
\ff{R}$ defined by the equation $z_2(\alpha)=f(z_1(\al)),$ has a
vertical tangent at the point $z(0)$. See the figure below for an example.
\begin{figure}
\centering
\includegraphics[width=.8\textwidth]{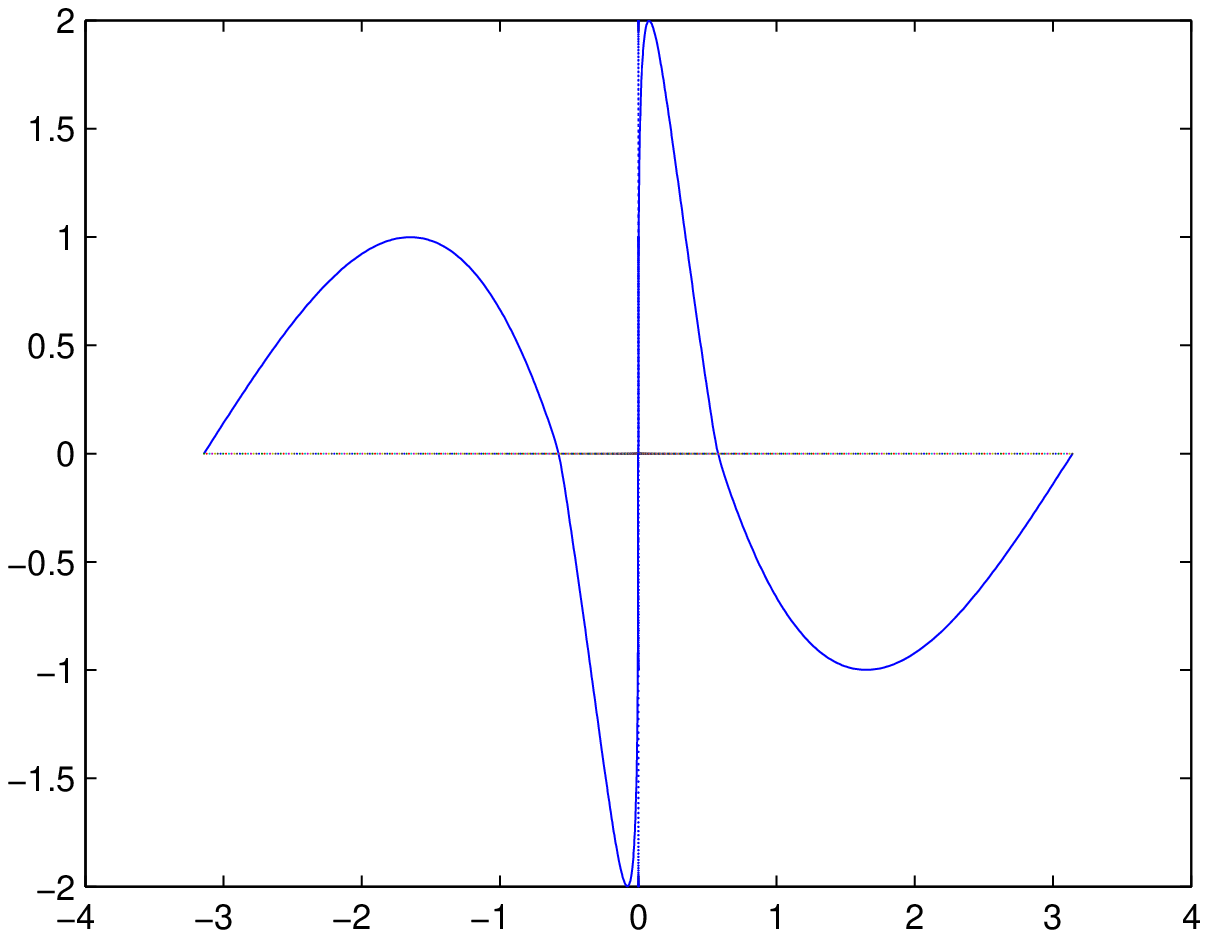}
\end{figure}
We shall prove the following
lemma:

\begin{lemma}\label{signo}
There exists a curve $z(\alpha)=(z_1(\alpha), z_2(\alpha))$ with
the following  properties:
\begin{enumerate}
\item $z_1(\al)-\alpha$ and $z_2(\al)$ are analytic
$2\pi-periodic$ functions and $z(\al)$ satisfies the arc-chord
condition, \item \label{odd}$z(\al)$ is  odd and \item
\label{vertical} $\da z_1(\alpha)>0$ if $\alpha \neq 0$, $\da
z_1(0)=0$ and $\da z_2(0)>0$,
\end{enumerate}

such that

\begin{eqnarray}\label{negatividad}
(\da v_1)(0)&=& \left.\left (\da
\int_{-\pi}^\pi\frac{\sin(z_1(\al)-z_1(\be))}{\cosh(z_2(\alpha)-z_2(\beta))-\cos(z_1(\alpha)-z_1(\beta))}(\da
z_1(\al)-\da z_1(\be))d\be \right)\right |_{\alpha=0} \nonumber \\&<& 0.
\end{eqnarray}
\end{lemma}

Proof:  We shall assume that $z(\alpha)$ is a
smooth curve satisfying the properties \ref{odd} and
\ref{vertical}. Differentiating the expression for the horizontal
component of the velocity, it is easy to obtain
$$(\da v_1)(\al)=\da\int_{-\pi}^\pi\frac{\sin(z_1(\al)-z_1(\al-\be))}
{\cosh(z_2(\alpha)-z_2(\alpha-\beta))-\cos(z_1(\alpha)-z_1(\alpha-\beta))}
(\da z_1(\al)-\da z_1(\al-\be))d\be$$
$$=\int_{-\pi}^\pi \frac{\cos(z_1(\alpha)-z_1(\alpha-\beta))(\da z_1(\al)-\da z_1(\al-\be))^2}
{\cosh(z_2(\alpha)-z_2(\alpha-\beta))-\cos(z_1(\alpha)-z_1(\alpha-\beta))}d\be$$
$$+\int_{-\pi}^\pi \frac{\sin(z_1(\al)-z_1(\al-\be))(\da^2 z_1(\al)-\da^2 z_1(\al-\be))}
{\cosh(z_2(\alpha)-z_2(\alpha-\beta))-\cos(z_1(\alpha)-z_1(\alpha-\beta))}d\be$$
 $$-\int_{-\pi}^\pi \sin((z_1(\al)-z_1(\al-\be)))(\da z_1(\al)-\da z_1(\al-\be))$$
 $$\times \frac{\sinh(z_2(\alpha)-z_2(\alpha-\beta))(\da z_2(\alpha)-\da z_2(\alpha-\beta))}
 {(\cosh(z_2(\alpha)-z_2(\alpha-\beta))-\cos(z_1(\alpha)-z_1(\alpha-\beta)))^2}d\be $$
$$-\int_{-\pi}^\pi \sin((z_1(\al)-z_1(\al-\be)))(\da z_1(\al)-\da z_1(\al-\be))$$
 $$\times \frac{\sin(z_2(\alpha)-z_2(\alpha-\beta))(\da z_1(\alpha)-\da z_1(\alpha-\beta))}
 {(\cosh(z_2(\alpha)-z_2(\alpha-\beta))-\cos(z_1(\alpha)-z_1(\alpha-\beta)))^2}d\be.$$
Evaluating this expression at $\al=0$ we have that
$$(\da v_1)(0)=\int_{-\pi}^\pi\frac{\cos(z_1(\be))(\da z_1(\be))^2+\sin(z_1(\be))\da^2 z_1(\be)}
{\cosh(z_2(\beta))-\cos(z_1(\beta))}d\be$$
$$-\int_{-\pi}^\pi \sin(z_1(\be))\da z_1(\be)\frac{\sin(z_1(\be))\da z_1(\be)-\sinh(z_2(\be))(\da z_2(0)-\da z_2(\be))}
{(\cosh(z_2(\beta))-\cos(z_1(\beta)))^2}d\be.$$

Integration by parts yields

$$\int_{-\pi}^\pi\frac{\sin(z_1(\be))\da^2
z_1(\be)}{\cosh(z_2(\beta))-\cos(z_1(\beta))}d\be$$
$$=-\int_{-\pi}^\pi\cos(z_1(\beta))\frac{(\da
z_1(\be))^2}{\cosh(z_2(\beta))-\cos(z_1(\beta))}d\be$$
$$+\int_{-\pi}^\pi \sin(z_1(\be))\da z_1(\be)\frac{\sin(z_1(\be))\da z_1(\be)+\sinh(z_2(\be))\da z_2(\be)}
{(\cosh(z_2(\beta))-\cos(z_1(\beta)))^2}d\be.$$

The above integrals converge because $z_1$ and $z_2$ satisfy the properties 2 and 3.
Therefore we obtain that
$$(\da v_1)(0)=\da
z_2(0)\int_{-\pi}^{\pi}\frac{\sin(z_1(\be))\sinh(z_2(\be))}{(\cosh(z_2(\beta))-\cos(z_1(\beta)))^2}\da
z_1(\be)d\be$$
\begin{equation}\label{reducida}
=2\da
z_2(0)\int_{0}^{\pi}\frac{\sin(z_1(\be))\sinh(z_2(\be))}{(\cosh(z_2(\beta))-\cos(z_1(\beta)))^2}\da
z_1(\be)d\be
\end{equation}

>From the expression (\ref{reducida}) we can control the sign of
$(\da v_1)(0)$. In order to clarify the proof we shall take
$$z_1(\be)=-\sin(\beta)+\beta.$$
We construct the function $z_2(\be)$ in the following way:

Let $\be_1$ and $\be_2$ be real numbers satisfying
$0<\be_1<\be_2<\pi$, and let $z^*(\be)$ be a smooth function on
$[-\pi, \pi]$, with the following properties,
$$
\begin{array}{ll}
a.\,\,z^*(\be)\mbox{ is odd.}\qquad & b.\,\,(\partial_\be z^*)(0)>0.\\
&\\
c.\,\,z^*(\be)>0 \mbox{ if } \be\in (0,\be_1). & d.\,\, z^*(\be)<0
\mbox{ if }\be\in (\be_1,\be_2]\\
&\\ e.\,\, z^*(\beta)\leq 0 \,\,\mbox{if}\,\, \beta\in
[\beta_2,\pi].
\end{array}
$$
For a positive real number $b$ to be fixed later, we define a
piecewise smooth function $\tilde{z}(\be)$ on $[-\pi,\pi]$, by
setting
\begin{eqnarray*}
\tilde{z}(\beta)&=& bz^*(\beta)\quad\text{if $|\be|\leq \be_1$},\\
\tilde{z}(\beta)&=& z^*(\beta) \quad\text{if $\be_1<|\beta|<\pi$}.
\end{eqnarray*}
Then
$$\int_{\be_1}^\pi\frac{\sin(z_1(\be))\sinh(\tilde{z}(\be))}{(\cosh(\tilde{z}(\beta))-\cos(z_1(\beta)))^2}\da
z_1(\be)d\be$$ is negative and independent of $b$, while
$$\int_0^{\be_1}\frac{\sin(z_1(\be))\sinh(\tilde{z}(\be))}{(\cosh(\tilde{z}(\beta))-\cos(z_1(\beta)))^2}\da
z_1(\be)d\be$$ tends to zero as $b\to\infty$.

Therefore, we can fix $b$ large enough so that
$$\int_0^{\pi}\frac{\sin(z_1(\be))\sinh(\tilde{z}(\be))}{(\cosh(\tilde{z}(\beta))-\cos(z_1(\beta)))^2}\da
z_1(\be)d\be<0.$$ It is now easy to approximate $\tilde{z}(\be)$
in $L^2[-\pi,\pi]$ by an odd, real-analytic $2\pi-$periodic
function such that
$$\int_0^{\pi}\frac{\sin(z_1(\be))\sinh(z_2(\be))}{(\cosh(z_2(\beta))-\cos(z_1(\beta)))^2}\da
z_1(\be)d\be<0,$$ and $\pa_\alpha z_2(0)>0$.

The conclusions of lemma (\ref{signo}) follow, thanks to
(\ref{reducida}).

 Theorem (\ref{kovalevsky}) and lemma
(\ref{signo}) allow us to show the breakdown of the R-T condition.

\begin{thm}
Let $z_0$ a curve satisfying the requirements of lemma
(\ref{signo}). Then there exists an analytic solution of the
Muskat problem satisfying the arc-chord condition in some interval
$[-T,T]$  such that for small enough  $T>0$ we have that:
\begin{enumerate} \item $\da z_1(\alpha,-t)>0$ $\forall \alpha$
and \item $\da z_1(0,t)<0$
\end{enumerate}
for all $t\in (0,T]$. In addition $\da z_2(0,t)>0$ in $[-T,T]$.
\end{thm}
Proof: We use theorem (\ref{kovalevsky}) to obtain the existence
and from lemma (\ref{signo}) we have that
$$(\pa_t\da z_1)(0,0)<0.$$
\begin{rem}
For $t\in[-T,0]$, our solution satisfies
$$\min_\alpha \pa_\alpha z_1(\alpha,t)>c|t|.$$
This follows easily, since $\pa_\alpha z_1(\alpha,0)$ has a
non-degenerate minimum at $\alpha=0$, and $\pa_t\pa_\alpha
z_1(0,0)<0$.
\end{rem}

\section{From a curve in $H^4$ in the stable regime to an analytic curve in the unstable regime}
Finally we show that there exists an open set of initial data in
the $H^4$ topology satisfying the arc-chord and R-T conditions
such that the solution for the Muskat problem reaches the unstable
regime. This section is devoted to proving theorem \eqref{oneone}.

Proof of theorem \eqref{oneone}: The idea is simply to take a small $H^4$-neighborhood of the initial data of an analytic solution.
Let $z_0$ be a curve as in lemma (\ref{signo}). Let
$z(\alpha,t)$ with $t\in [-T,T]$ for some $T>0$, the solution for
the equation (\ref{ec1d}) given by theorem \ref{kovalevsky}. We
consider the curve
$w^\ep_\delta(\alpha)=((w^\ep_\delta(\alpha)_1,(w^\ep_\delta(\alpha))_2)$
which is a small perturbation in $H^4(\T)$ of the curve
$z(\alpha,t)$ at time $t=-\delta$, with $0<\delta<T$, i. e.
$$||w_\delta^\ep(\cdot)-z(\cdot,-\delta)||_{H^4}=||\eta_\delta^\ep||_{H^4}\leq
\ep.$$ Also, $w_\delta^\ep(\alpha)$ satisfies the R-T condition
$$\sigma_{w_\delta^\ep}(\alpha)\equiv(\rho^2-\rho^1)\da
(w_\delta^\ep)_1(\alpha)> 0,$$ if $0<\delta\leq \delta_0$ and
$0<\ep\leq \ep(\delta)$. From now on, we take $\ep$ and $\delta$
to satisfy this condition. Also, we may take
$\ep(\delta)<\ep_0$.

Since $z(\alpha,0)=z_0(\alpha)$ is a smooth curve satisfying the
arc-chord condition we can assume that there exist $\ep_0>0$ and
$0<\delta_0<T$ such that
\begin{equation}\label{cota1}\sup_{0<\ep\leq \ep_0,\,0<\delta\leq
\delta_0} ||w_\delta^\ep(\cdot)||_{H^4(\T)}\leq
C(z_0,\ep_0,\delta_0),\end{equation} and
\begin{equation}\label{cota2}\sup_{0<\ep\leq \ep_0,\,0<\delta\leq \delta_0}
||F(w_\delta^\ep)||_{L^\infty(\T)}\leq
C(z_0,\ep_0,\delta_0).\end{equation}

Now, let the curve $w^\ep(\alpha,t)$ be the solution to the
equation
\begin{align*}
\partial_t w^\ep (\al,t) = &\int
\frac{\sin((w^\ep)_1(\al,t)-(w^\ep)_1(\beta,t))}
{\cosh((w^\ep)_2(\al,t)-(w^\ep)_2(\beta,t))-\cos((w^\ep)_1(\al,t)-(w^\ep)_1(\beta,t))}\\
& \times (\da
w^\ep(\al,t) - \da w^\ep(\beta,t)) d\beta\\
w^\ep(\alpha,-\delta)&= w_{\delta}^\ep(\alpha).
\end{align*}
>From theorems \ref{instantanalycity} and \ref{RTmoe0} and the
inequalities (\ref{cota1}) and (\ref{cota2}) we see that we can
choose $\ep_0$ and $\delta$ small enough in such a way that
$w^\ep(\alpha,t)$ is well defined, for all $0<\ep<\ep_0$, in
$t\in[-\delta,0]$ unless $w^\ep(\alpha,t)$ loses the R-T
condition. That means  there exist some point $\alpha_0$ and some
time $t_0\in [-\delta,0]$ satisfying
$\sigma_{w^\ep}(\alpha_0,t_0)<0.$ Also, for small enough $\ep_0$
and fixed $\delta$ we have that
\begin{equation}\label{rai}
\sigma_{w^\ep}(\alpha)\geq a>0,
\end{equation}
where $a$ is a real number independent of $\ep$. The numbers
$\ep_0$ and $\delta$ are fixed for the rest of the proof.

If there exist times $t$ such that there exists some point
$\alpha_0$ with $\sigma_{w^\ep}(\alpha_0,t)=0$, we denote
the first of these times to be $T^\ep\in (-\delta,\infty)$. Also
we set  $\tilde{T}^\ep=\min \{T^\ep,0\}$ and
$I^\ep=[-\delta,\tilde{T}^\ep]$. Due to (\ref{rai}) we have that
$$\inf_{0<\ep<\ep(\delta)}\tilde{T}^\ep>t_b>-\delta,$$
for some  number $t_b$.

>From the proof of theorems \ref{instantanalycity} and \ref{RTmoe0}
we know that there exists a function $h(t)$, given by the
expression
\begin{equation*}h(t)=\left\{\begin{array}{ccc}c_a(t+\delta) &
-\delta\leq t \leq t_a\\
c_a (t_a+\delta) e^{-C_a(t-t_a)} &
t>t_a,\end{array}\right.\end{equation*} where $t_a$ (small enough),
$c_a$ and $C_a$ are constants which only depend on the constant
$C(z_0,\ep_0,\delta)$ (see (\ref{cota1}) and (\ref{cota2})), such
that $w^\ep(\alpha,t)$ is  an analytic function in the strip
$$S(t)=\{ \zeta\in \C\,:\, |\Im(\zeta)|< h(t)\},$$ and
also
$$(||w^\ep(\cdot,t)||_S+1)^k\leq C_a,$$
for some large enough $k$ and $t\in[t_a,\tilde{T}^\ep]$ (notice
that the constants $t_a$, $c_a$ and $C_a$ do not depend on $\ep$).

In this situation we claim the following:
$$
\frac{d}{dt}\int_{-\pi}^{\pi}|\pa_\alpha^4(w^\ep(\alpha\pm i
h(t),t)-z(\alpha\pm i h(t),t))|^2d\alpha$$
\begin{equation}\label{diferencia}\leq
C(||\pa_\alpha^4(w^\ep(\cdot\pm i h(t),t)-z(\cdot\pm i
h(t),t))||^2_{L^2(\T)}+||w^\ep(\cdot+i h(t),t)-z(\cdot+ i
h(t),t)||_{L^2(\T)}^2)\end{equation} for $t\in I^\ep$ and where $C$
is a constant just depending on $C(z_0,\ep_0,\delta)$.

We will prove this inequality at the end of the section. Let us
assume that (\ref{diferencia}) holds.

We notice that we  can always choose either  a subsequence
$\{\ep_n\}_1^\infty$ with  $\ep_n\to 0$ when $n\to \infty$ such
that $T^{\ep_n}<0$ $\forall n$ or a subsequence
$\{\ep_m\}_1^\infty$ with $\ep_m\to 0$ when $m\to \infty$ such
that $T^{\ep_m}\geq 0$ $\forall m$ (the case in which there exist
only a finite number of times $T^\ep$ can be treated as this last
case). We deal with these two cases, I and II, separately:

I. $T^\ep<0$ for all  $\ep$. From inequality (\ref{diferencia}) we
can take $\ep$ small enough such that
$$w^\ep(\alpha,T^\ep)- z(\alpha,T^\ep)$$
has norm $\leq C\ep$ in $H^4(S(T^\ep))$.

Note that
$$0=\min_{\alpha} \pa_\alpha (w^\ep)_1(\alpha,T^\ep)\geq
-C\ep+\min_\alpha \pa_\alpha z^0_1\geq -C\ep+c|T^\ep|$$ by the
remark at the end of section 5. Thus, $|T^\ep|<C\ep$.

 Then
$$z(\alpha,T^\ep)- z_0(\alpha)$$
has norm $\leq C\ep$ in $H^4(S(0))$ and therefore
$$|(\da (v_{w^\ep})_1)(\alpha_0, T^{\ep})- (\da (v_{z_0})_1)(0)|\leq C\ep,$$
 and we can conclude that
$$(\da (v_{w^\ep})_1)(\alpha_0, T^{\ep})<0.$$

Here we recall that
$$(v_{w^\ep})_1=\partial_t w^\ep (\al,t) =\int
\frac{\sin((w^\ep)_1(\al,t)-(w^\ep)_1(\beta,t))(\da w^\ep(\al,t) - \da w^\ep(\beta,t)) d\beta}{\cosh((w^\ep)_2(\al,t)
-(w^\ep)_2(\beta,t))-\cos((w^\ep)_1(\al,t)-(w^\ep)_1(\beta,t))}$$

Applying
the same argument as in  section \ref{AsAu} to the curve
$w^\ep(\alpha, T^\ep)$ we finish the proof of theorem
\ref{oneone} in the case $T^\ep<0$ for all $\ep$.

II. $T^\ep\geq 0$ for all  $\ep$. Then we can apply a
Cauchy-Kowalewski theorem to the initial data
$$w^\ep(\alpha,0)-z(\alpha,0)$$ satisfying
$$||w^\ep-z||_{S(0)}\leq
C\ep.$$  For $t>0$ small enough, $z(\alpha,t)$ is in the unstable
regime. We achieve the conclusion  of theorem \ref{oneone}
by continuity with respect to the initial data.

The rest of the section is devoted to proving inequality
(\ref{diferencia}). We shall denote $\gamma=\alpha + i h(t)$ and
$d(\gamma,t)=\da^4 (w(\gamma,t)-z(\gamma,t)$) (we omit the
superscript $\ep$ in the notation) and we recall that $w(\al,t)$
and $z(\alpha,t)$ are real for real $\alpha$ (therefore we obtain
similar similar estimates for $\gamma=\alpha-i h(t)$). In order to
prove inequality (\ref{diferencia}) we have to compute the
following quantity
\begin{align*}
&\frac{d}{dt}\int_{-\pi}^{\pi}|d(\gamma,t)|^2d\al&= 2\Re \left\{\int_{-\pi}^\pi
\overline{d(\gamma,t)}d_t(\gamma,t)d\alpha\right\}
+2\Re\left\{ih'(t)\int_{-\pi}^{\pi}\overline{d(\gamma,t)}\pa_\alpha
d(\gamma,t)d\al\right\}.
\end{align*}
Again we treat in detail the most singular term in
$d_t(\gamma,t)$. Recall  $K(\alpha,\beta)$  from section 3 and
write $K_w$ and $K_z$ for corresponding expressions arising from
$z$ and $w$. Then we have that
\begin{align*}
d_t(\gamma,t)=&
\int_{-\pi}^{\pi}K_w(\gamma,\gamma-\beta)\da^5(w(\gamma,t)-w(\gamma-\beta,t)d\be\\
&-\int_{-\pi}^{\pi}K_z(\gamma,\gamma-\beta)\da^5(z(\gamma,t)-z(\gamma-\beta,t))d\be+\text{l.o.t}(\alpha,t),
\end{align*}
where
$$2\Re \left\{\int_{-\pi}^\pi
\overline{d(\gamma,t)}\text{l.o.t}(\alpha)d\alpha\right\}\leq
C(||d(\cdot+ih(t),t)||_{L^2(\T)}^2+||w(\cdot+i h(t),t)-z(\cdot+ i
h(t),t)||_{L^2(\T)}^2).$$ Here $C$ is a constant which just
depends on $\ep_0$ and $\delta$.

We can write
\begin{align*}
&
\int_{-\pi}^{\pi}K_w(\gamma,\gamma-\beta)\da^5(w(\gamma,t)-w(\gamma-\beta,t))d\be\\
&\qquad\qquad\qquad\qquad-\int_{-\pi}^{\pi}K_z(\gamma,\gamma-\beta)\da^5(z(\gamma,t)-z(\gamma-\beta,t))d\be
\end{align*}
\begin{align*}
=&\int_{-\pi}^{\pi}K_w(\gamma,\gamma-\beta)\da^5((w(\gamma,t)-z(\gamma,t))-(w(\gamma-\beta,t)-z(\gamma-\beta,t)))d\be\\
&+\int_{-\pi}^\pi
\left\{K_z(\gamma,\gamma-\beta)-K_w(\gamma,\gamma-\beta)\right\}\da^5(z(\gamma,t)-z(\gamma-\beta,t))d\beta
\end{align*}
\begin{align*}
=&\int_{-\pi}^{\pi}K_w(\gamma,\gamma-\beta)\pa_\alpha
(d(\gamma,t)-d(\gamma-\beta,t)) d\beta\\
&+\int_{-\pi}^\pi
\left\{K_z(\gamma,\gamma-\beta)-K_w(\gamma,\gamma-\beta)\right\}\da^5(z(\gamma,t)-z(\gamma-\beta,t))d\beta\\
\equiv& X_1(\alpha,t)+X_2(\alpha,t).
\end{align*}
Therefore
\begin{align*}
\frac{d}{dt}\int_{-\pi}^{\pi}|d(\gamma,t)|^2d\al\leq & C||d(\cdot+ih(t),t)||_{L^2(\T)}^2 \\
&+2\Re \left\{\int_{-\pi}^\pi
\overline{d(\gamma,t)}X_1(\alpha,t)d\alpha\right\}+ 2\Re
\left\{\int_{-\pi}^\pi
\overline{d(\gamma,t)}X_2(\alpha,t)d\alpha\right\}\\
&+2\Re\left\{ih'(t)\int_{-\pi}^{\pi}\overline{d(\gamma,t)}\pa_\alpha
d(\gamma,t)d\al\right\}.
\end{align*}
Following the computations in section 3 when $t\in[-\delta,t_a]$
and those in section 4 when $t\in[t_a,\tilde{T}^\ep]$ we have that
\begin{align*}
&\frac{d}{dt}\int_{-\pi}^{\pi}|d(\gamma,t)|^2d\al\leq
C||d(\cdot+ih(t),t)||_{L^2(\T)}^2+ 2\Re \left\{\int_{-\pi}^\pi
\overline{d(\gamma,t)}X_2(\alpha,t)d\alpha\right\}.
\end{align*}
In addition
\begin{align*}
&\left|\frac{\sin(w^1(\gamma)-w^1(\gamma-\beta))}
{\cosh(w^2(\gamma)-w^2(\gamma-\beta))-\cos(w^1(\gamma)-w^1(\gamma-\beta))}\right.\\
-&\left.\frac{\sin(z^1(\gamma)-z^1(\gamma-\beta))}
{\cosh(z^2(\gamma)-z^2(\gamma-\beta))-\cos(z^1(\gamma)-z^1(\gamma-\beta))}\right|\\
=&\left|\left\{\left [K_w(\gamma,\gamma-\beta)-\frac{\pa_\al
w^1(\gamma)}{(\pa_\al w^1(\gamma))^2+(\pa_\al
w^1(\gamma))^2}\cot\left(\frac{\be}{2}\right)\right]\right.\right.\\
-& \left.\left[K_z(\gamma,\gamma-\beta)-\frac{\pa_\al
z^1(\gamma)}{(\pa_\al z^1(\gamma))^2+(\pa_\al
z^1(\gamma))^2}\cot\left(\frac{\be}{2}\right)\right]\right\}\\
+ & \left.\left\{\frac{\pa_\al w^1(\gamma)}{(\pa_\al
w^1(\gamma))^2+(\pa_\al w^1(\gamma))^2}-\frac{\pa_\al
z^1(\gamma)}{(\pa_\al z^1(\gamma))^2+(\pa_\al
z^1(\gamma))^2}\right\}\cot\left(\frac{\be}{2}\right)\right|\\
\leq &
(||d(\cdot+ih(t),t)||_{L^2(\T)}+||w(\cdot+i h(t),t)-z(\cdot+ i
h(t),t)||_{L^2(\T)})\left\{C+C\left|\cot\left(\frac{\be}{2}\right)\right|\right\}.
\end{align*}
Also, $$|\pa_\alpha^5 z(\alpha\pm i h(t),t)-\pa_\alpha^5
z(\alpha\pm i h(t)-\beta,t)|\leq C
\left|\tan\left(\frac{\be}{2}\right)\right|$$ since $z$ is the
analytic unperturbed solution. Therefore
$$2\Re\left\{\int_{-\pi}^\pi
\overline{d(\gamma,t)}X_2(\alpha,t)d\alpha\right\}\leq
C(||d(\cdot+ih(t),t)||_{L^2(\T)}^2+||w(\cdot+i h(t),t)-z(\cdot+ i
h(t),t)||_{L^2(\T)}^2).$$ We are done.

\section{Turning water waves}\label{ww}

 Let us consider an incompressible irrotational flow satisfying the Euler equations
\begin{equation}\label{Euler}
\rho(v_t+v\cdot\grad v)=-\grad p-g\rho(0,1),
\end{equation}
where $\rho$ satisfies (\ref{trasporte},\ref{frho}) and $\rho^1=0$.
This system of equations provides the motion of the interface for the
water wave problem (see \cite{BL, Lannes} and references therein),
whose contour equation is given by
\begin{equation}\label{vww}
z_t(\al,t)=BR(z,\omega)(\al,t)+c(\al,t)\da z(\al,t),
\end{equation}
and
\begin{align}
\begin{split}\label{cEuler}
\omega_t(\al,t)&=-2\partial_t BR(z,\omega)(\al,t)\cdot
\da z(\al,t)-\da( \frac{|\omega|^2}{4|\da z|^2})(\al,t) +\da (c\, \omega)(\al,t)\\
&\quad +2c(\al,t)\da BR(z,\omega)(\al,t)\cdot\da z(\al,t)-2g\da
z_2(\al,t).
\end{split}
\end{align}
The values of $z(\alpha,t)$ and $w(\alpha,t)$ are given at an
initial time $t_0$: $z(\alpha,t_0)=z^0(\alpha)$ and
$w(\alpha,t_0)=w^0(\alpha)$. For more details see \cite{ADP2}.

As an application of section 5, we can consider initial data given
by a graph $(\al,f_0(\al))$ and show that in finite time the
interface evolution reaches a regime where the contour only can be
parametrized as $z(\al,t)=(z_1(\al,t),z_2(\al,t)),$ for
$\al\in\R$, with $\da z_1(\al,t)<0$ for  $\al\in I$, a non-empty
interval. This implies that there exists a time $t^*$ where the
solution of the free boundary problem reparametrized by
$(\al,f(\al,t))$ satisfies $\|f_\al\|_{L^\infty}(t^*)=\infty$.
\begin{thm}\label{waterwaves}
There exists a non-empty open set of initial data
$z^0(\alpha)=(\al,f_0(\al))$ and $w^0(\alpha)$, with $f_0\in H^5$
and $w^0\in H^4$, such that in finite time $t^*$ the solution of
the water wave problem (\ref{vww},\ref{cEuler}) given by
$(\al,f(\al,t))$ satisfies $\|f_\al\|_{L^\infty}(t^*)=\infty$. The
solution can be continued for $t>t^*$ as $z(\al,t)$ with $\da
z_1(\al,t)<0$ for $\al\in I$, a non-empty interval.
\end{thm}

Proof: Let us consider a curve $z^*(\al)\in H^5$ satisfying 1., 2.
and 3. of Lemma \ref{signo}. We point out that  analyticity  is not required here. In order to find a velocity with
property \eqref{negatividad} we pick for water waves
$\omega(\al,t^*)=-\da z^*_2(\al)$ and a suitable
$z(\al,t^*)=z^*(\al)$ as an initial datum. Notice that the
tangential term does not affect the evolution. Then, with the
appropriate $c(\al,t)$, we can apply the local existence result in
\cite{ADP2}: There exists a solution of the water wave problem
with $z(\al,t)\in C([t^*-\delta,t^*+\delta];H^5)$,
$\omega(\al,t)\in C([t^*-\delta,t^*+\delta];H^4)$ and $\delta>0$
small enough. The initial data promised by theorem \ref{waterwaves} are any
sufficiently small perturbations of $z(\alpha,t)$ and $w(\alpha,t)$
at time $t=t^*-\delta$.

\section{Breakdown of Smoothness}

In \cite{ADCPM3} we will exhibit a solution $z(\al,t)$ of the Muskat equation, with the following properties.
\begin{enumerate}
\item At time $t_0$, the interface is real-analytic and satisfies the arc-chord and Rayleigh-Taylor conditions.
\item At time $t_1>t_0$, the interface turns over.
\item At time $t_2>t_1$, the interface no longer belongs to $C^4$, although it is real-analytic for all times $t\in[t_0,t_2)$.
\end{enumerate}
In this section we provide a brief sketch of our proof of the existence of such a Muskat solution.

Our Muskat solution $z(\al,t)$ will be a small perturbation of a Muskat solution $z^{00}(\al,t)$, with the following properties.
\begin{itemize}
\item [4.] $z^{00}(\alpha,t)$ is real analytic in $\alpha$, for $|\Im \alpha|<\ep^{00}$ and $|\tau|\leq \tau^{00}$.
\item [5.] For $t\in [-\tau^{00},0)$, $z^{00}(\alpha,t)$ satisfies the Rayleigh-Taylor and arc-chord conditions.
\item [6.] For $t=0$, the curve $z^{00}(\alpha,t)$ has a vertical tangent at $\alpha=0$.
\item [7.] For $t\in (0,\tau^{00}]$, the curve $z^{00}(\alpha,t)$ fails to satisfy the Rayleigh-Taylor condition.
\end{itemize}
This paper constructs Muskat solutions $z^{00}$ satisfying 4., 5., 6. and 7. Our problem is to pass from $z^{00}$ to a nearby Muskat solution $z$ satisfying 1., 2. and 3. The idea is as follows.

So far, we have studied the analytic continuation of Muskat solutions to a time-varying strip
\begin{equation}\nonumber
S(t)=\{|\Im \alpha|\leq h(t)\},
\end{equation}
in the complex plane. In our forthcoming paper \cite{ADCPM3}, we will study the analytic continuation of a Muskat solution to a carefully chosen time-varying domain of the form
\begin{equation}\label{8.8}
\Omega(t)=\{|\Im \alpha|\leq h(\Re\al,t)\},
\end{equation}
defined for $t\in[-\tau^{10},\tau].$ Here, $\tau$ is a small enough positive number.

For $t\in [-\tau^{10},\tau]$, we will work with the space $H^4(\Omega(t))$, consisting of all analytic functions $F\,:\,\Omega(t)\mapsto \C^2$ whose derivatives up to order $4$ belong to $L^2(\pa\Omega(t))$.

We will pick our time-varying domain $\Omega(t)$ in \eqref{8.8} so that $h(x,t)>0$ for all $(x,t)\in \R/2\pi\Z\times[-\tau^{10},\tau)$ and $h(x,\tau)>0$ for all $x\in \R/2\pi\Z\setminus \{0\}$, but $h(0,\tau)=0$. Thus, the domain $\Omega(t)$ has 'thickness' zero at the origin. Consequently, $H^4(\Omega(\tau))$ is not contained in $C^4(\R/2\pi\Z)$.

We will also take $\tau<\tau^{00}$ and $h(x,t)<\ep^{00}$, so that the Muskat solution $z^{00}(\alpha,t)$ continues analytically to $\Omega(t)$,  for each $t\in[-\tau^{10},\tau]$.

We can therefore pick an 'initial' curve $z^0(\alpha)$, such that
\begin{itemize}
\item[8.] $z^0(\alpha)-z^{00}(\alpha,\tau)$ belongs to $H^4(\Omega(\tau))$ and has small norm, yet
\item[9.] $z^0(\al)$ does not belong to $C^4(\R/2\pi\Z)$.
\end{itemize}
We solve the Muskat problem backwards in time, with the 'initial' condition
\begin{itemize}
\item[10.] $z(\alpha,\tau)=z^0(\alpha)$.
\end{itemize}
By a more elaborate version of the analytic continuation arguments used in this paper, we find that our Muskat solution exists and continues analytically into $\Omega(t)$, for all $t\in[t_*,\tau]$ (for a suitable time $t_*$); moreover,
\begin{itemize}
\item[11.] $z(\al,t)-z^{00}(\alpha,t)$ has small norm in $H^4(\Omega(\tau))$, for all $t\in[t_*,\tau]$.
\end{itemize}
Here, either
\begin{itemize}
\item[12.]$t_*=-\tau^{10}$ or
\item[13.] a modified Rayleigh-Taylor condition, adapted to the time-varying domain, fails at time $t_*$.
\end{itemize}
We can rule out 13., thanks to 11., together with our understanding of $z^{00}(t)$ and $\Omega(t)$.

Thus, we obtain a Muskat solution $z(\alpha,t)$, satisfying 9., 10., 11. and 12. Properties 1., 2. and 3. of $z(\alpha,t)$ now follow easily.

\subsection*{{\bf Acknowledgements}}

\smallskip

 AC, DC and FG were partially supported by the grant {\sc MTM2008-03754} of the MCINN (Spain) and
the grant StG-203138CDSIF  of the ERC. CF was partially supported by
NSF grant DMS-0901040 and ONR grant ONR00014-08-1-0678.
FG was partially supported by NSF grant DMS-0901810. MLF was partially supported by the grants {\sc MTM2008-03541} and {\sc MTM2010-19510} of the MCINN (Spain). The authors thank Professor Rafael de la Llave for helpful discussions.

\begin{tabular}{ll}
\textbf{Angel Castro} &  \\
{\small Instituto de Ciencias Matem\'aticas} & \\
{\small Consejo Superior de Investigaciones Cient\'ificas} &\\
{\small Serrano 123, 28006 Madrid, Spain} & \\
{\small Email: angel\underline{  }castro@icmat.es} & \\
   & \\
\textbf{Diego C\'ordoba} &  \textbf{Charles Fefferman}\\
{\small Instituto de Ciencias Matem\'aticas} & {\small Department of Mathematics}\\
{\small Consejo Superior de Investigaciones Cient\'ificas} & {\small Princeton University}\\
{\small Serrano 123, 28006 Madrid, Spain} & {\small 1102 Fine Hall, Washington Rd, }\\
{\small Email: dcg@icmat.es} & {\small Princeton, NJ 08544, USA}\\
 & {\small Email: cf@math.princeton.edu}\\
 & \\
\textbf{Francisco Gancedo} &  \textbf{Mar\'ia L\'opez-Fern\'andez}\\
{\small Department of Mathematics} & {\small Institut f\"ur Mathematik}\\
{\small University of Chicago} & {\small Universit\"at Z\"urich}\\
{\small 5734 University Avenue,} & {\small Winterthurerstr. 190,
CH-8057}\\
{\small Chicago, IL 60637, USA}  & {\small Z\"urich, Switzerland}\\
{\small Email: fgancedo@math.uchicago.edu} & {\small Email:
maria.lopez@math.uzh.ch}
\end{tabular}


\begin{thebibliography}{99}
\bibitem{Am} D. Ambrose. Well-posedness of Two-phase Hele-Shaw Flow without Surface Tension.
\emph{Euro. Jnl. of Applied Mathematics} 15, (2004), 597-607.

\bibitem{Bardos} C. Bardos and S. Benachour. Domaine d'analyticite des solutions
de l'equation d'Euler dans un ouvert de $\R^n$.
\emph{Annal. Sc. Normale Sup. di Pisa}, (1977), no. 4, 647-687.

\bibitem{BL} C. Bardos and D. Lannes.  Mathematics for 2d Interfaces.
ArXiv:1005.5329.  \emph{To appear in Panorama et Syntheses} (2010).

\bibitem{hou}J. Beale, T. Y. Hou and J. Lowengrub. Convergence of a boundary integral method for water waves.
SIAM J. Numer. Anal. 33 , no. 5, 1797�1843 (1996).

\bibitem{ADCPM3} A. Castro, D. C\'ordoba, C. Fefferman and F. Gancedo. Breakdown of smoothness for the Muskat
problem. \emph{Preprint} (2011).

\bibitem{turning}A. Castro, D. C\'ordoba, C. Fefferman, F. Gancedo and Mar\'ia L\'opez-Fern\'andez.
Turnning waves and breakdown for incompressible flows. \emph{Proc. Natl. Acad. Sci.} 108, no. 12, 4754-4759 (2011).

\bibitem{CE} A. Constantin and J. Escher. Wave breaking for nonlinear nonlocal shallow water equations. \emph{Acta Math.},
181, 229-243 (1998).

\bibitem{ccgs} P. Constantin, D. C\'ordoba, F. Gancedo and R.M. Strain. On the global existence  for
the Muskat problem. \emph{To appear in JEMS} (2011).

\bibitem{Peter} P. Constantin and M. Pugh. Global solutions for small data to the
Hele-Shaw problem. \emph{Nonlinearity}, 6 (1993), 393 - 415.

\bibitem{AD} A. C\'ordoba, D. C\'ordoba. A pointwise estimate for fractionary derivatives with applications
to partial differential equations. \emph{Proc. Natl. Acad. Sci. USA} 100, 26,
(2003), 15316-15317.

\bibitem{ADY} A. C\'ordoba, D. C\'ordoba and F. Gancedo. Interface evolution: the Hele-Shaw and Muskat problems.
 \emph{Annals of Math.}, 173, 1, (2011), 477-542.

\bibitem{ADP2} A. C\'ordoba, D. C\'ordoba, F. Gancedo.
Interface evolution: the water wave problem in 2D. \emph{Adv. Math.}, 223,
no. 1, (2009), 120-173.

\bibitem{DY} D. C\'ordoba and F. Gancedo. Contour dynamics of incompressible 3-D fluids
in a porous medium with different densities. \emph{Comm. Math.
Phys.} 273, 2,(2007), 445-471.

\bibitem{DP2} D. C\'ordoba and F. Gancedo. A maximum principle for the Muskat problem for fluids with different densities.
\emph{Comm. Math.Phys.}, 286 (2009), no. 2, 681-696.

\bibitem{DP3} D. C\'ordoba and F. Gancedo. Absence of squirt singularities for the multi-phase Muskat problem.
\emph{Comm. Math. Phys.}, 299, 2, (2010), 561-575.

\bibitem{DPR} D. C\'ordoba, F. Gancedo and R. Orive. A note on the interface dynamics for convection in porous media.
\emph{Physica D}, 237 (2008), 1488-1497.

\bibitem{Drit} D.G. Dritschel. Contour dynamics and contour surgery: numerical
    algorithms for extended, high-resolution modelling of vortex dynamics in
    two-dimensional, inviscid, incompressible flows. \emph{Computer Physics Report}, 10, (1989),77-146.

\bibitem{Esch1} J. Escher and G. Simonett. Classical solutions for Hele-Shaw models with surface tension.
\emph{Adv. Differential Equations}, 2, (1997), 619-642.

\bibitem{Esch2} J. Escher and B.-V. Matioc. On the parabolicity of the Muskat problem: Well-posedness, fingering,
and stability results. \emph{Z. Anal. Awend.} 30, 193-218, (2011).

\bibitem{Esch3} J. Escher, A.V. Matioc and B.-V. Matioc. A generalised Rayleigh-Taylor condition for the Muskat problem.
ArXiv:1005.2511. \emph{Preprint} (2010).

\bibitem{Hai} E. Hairer and G. Wanner. Solving ordinary differential equations.
    I. Nonstiff problems. Springer Series in Computational Mathematics, 8.
    Springer-Verlag, Berlin, 1987.

\bibitem{H-S} Hele-Shaw. On the motion of a viscous fluid between two parallel plates. \emph{Trans. Royal
Inst. Nav. Archit.}, London 40, 21, (1898).

\bibitem{Howison} S. Howison. A note on the two-phase Hele-Shaw
problem. \emph{J. Fluid Mech.}, vol. 409, (2000), 243-249.

\bibitem{HLS} T. Hou, J. Lowengrubb and M. Shelley. Removing the stiffness
from interfacial flows with surface tension. \emph{J. Comput.
Phys.}, 114, (1994), 312-338.

\bibitem{Lannes} D. Lannes. A stability criterion for two-fluid interfaces and
applications. ArXiv:1005.4565. \emph{Preprint} (2010).

\bibitem{Muskat} M. Muskat. Two fluid systems in porous media. The encroachment of water into an oil sand.
\emph{Physics}, 5, (1934), 250-264.

\bibitem{Nirenberg} L. Nirenberg. An abstract form of the nonlinear Cauchy-Kowalewski theorem.
\emph{J. Differential Geometry}, 6, (1972), 561-576.

\bibitem{Nishida} T. Nishida. A note on a theorem of Nirenberg. \emph{J. Differential Geometry}, 12, (1977), 629-633.

\bibitem{Otto} F. Otto. Viscous fingering: an optimal bound on the growth rate of the mixing zone.
\emph{SIAM J. Appl. Math.}  57,  no. 4, (1997), 982-990.

\bibitem{Ray} Lord Rayleigh (J.W. Strutt), On the instability of jets. \emph{Proc. Lond. Math.
Soc.} 10, 413, (1879).

\bibitem{S-T} P.G. Saffman and Taylor. The penetration of a fluid into a porous medium or
Hele-Shaw cell containing a more viscous liquid.
\emph{Proc. R. Soc. London, Ser. A} 245, (1958), 312-329.

\bibitem{SCH} M. Siegel, R. Caflisch and S. Howison. Global
Existence, Singular Solutions, and Ill-Posedness for the Muskat
Problem. \emph{Comm. Pure and Appl. Math.}, 57, (2004), 1374-1411.

\bibitem{SSBF} C. Sulem, P.L. Sulem, C. Bardos and U. Frisch. Finite time analyticity for the two- and three-dimensional Kelvin-Helmholtz instability.
\emph{Comm. Math. Phys.} 80, 4, (1981), 485-516.

\bibitem{Yi} F. Yi. Local classical solution of Muskat free boundary problem, \emph{J. Partial Diff. Eqs.}, 9 (1996), 84-96.

\bibitem{Yi2} F. Yi. Global classical solution of Muskat free boundary problem, \emph{J. Math. Anal. Appl.}, 288 (2003), 442-461.

\end{thebibliography}
\end{document}